\documentclass[12pt,english]{amsart}
\usepackage[T1]{fontenc}
\usepackage[latin9]{inputenc}
\usepackage{geometry}
\usepackage{color}
\usepackage{babel}
\usepackage{verbatim}
\usepackage{amsthm}
\usepackage{amstext}
\usepackage{amssymb}
\usepackage[unicode=true]
 {hyperref}
\usepackage{breakurl}

\makeatletter
\numberwithin{equation}{section}
\numberwithin{figure}{section}
\theoremstyle{plain}
\newtheorem{thm}{\protect\theoremname}[section]
  \theoremstyle{definition}
  \newtheorem{defn}[thm]{\protect\definitionname}
  \theoremstyle{remark}
  \newtheorem{rem}[thm]{\protect\remarkname}
  \theoremstyle{plain}
  \newtheorem{fact}[thm]{\protect\factname}
  \theoremstyle{plain}
  \newtheorem{lem}[thm]{\protect\lemmaname}
  \theoremstyle{plain}
  \newtheorem{cor}[thm]{\protect\corollaryname}
  \theoremstyle{remark}
  \newtheorem*{rem*}{\protect\remarkname}

\usepackage{amsfonts}

\usepackage{amscd}

\makeatother

  \providecommand{\corollaryname}{Corollary}
  \providecommand{\definitionname}{Definition}
  \providecommand{\factname}{Fact}
  \providecommand{\lemmaname}{Lemma}
  \providecommand{\remarkname}{Remark}
\providecommand{\theoremname}{Theorem}

\begin{document}

\title[Duality and interpolation spaces]{Lecture notes on duality and interpolation spaces
\\
\scriptsize{(Version 2, Nov. 1, 2014)}
}

\author{Michael Cwikel}

\address{Department of Mathematics, Technion - Israel Institute of Technology,
Haifa 32000, Israel}

\email{mcwikel@{math.technion.ac.il}}

\thanks{This research was supported by the Technion V.P.R.\ Fund and by
the Fund for Promotion of Research at the Technion. I am also very
grateful for the hospitality of the Centre for Mathematics and its
Applications at the Australian National University.}

\subjclass{Primary 46B70, Secondary 46B10, } %% 46B70 = Interpolation between normed linear spaces.
%% 46B10 = Duality and reflexivity.

\keywords{Dual space, interpolation space, complex interpolation}

\begin{comment}
\textbf{\textcolor{red}{VERSION 1 OF THIS PAPER WAS POSTED ON THE
ARXIV ON 25/3/2008. }}
\end{comment}

\begin{abstract}
Known or essentially known results about duals of interpolation spaces
are presented, taking a point of view sometimes slightly different
from the usual one. Particular emphasis is placed on Alberto Calderón's
theorem describing the duals of complex interpolation spaces. The
pace is slow, since these notes are intended for graduate students
who have just begun to study interpolation spaces. 
\end{abstract}

\maketitle

This paper, or set of lecture notes, gives a careful, pedantic, perhaps
even a bit too pedantic, treatment of the general topic of duality
in the context of Banach couples. I will follow the same standard
(it might be called ``naive'') approach which is explicitly or implicitly
used in many papers about interpolation spaces, including the classical
expositions of Alberto Calderón \cite{ca} and Lions--Peetre \cite{lp}.
It should be mentioned that 
Sten Kaijser and Joan Wick--Pelletier (see \cite{kp1}, \cite{kp2},
\cite{kp3} and \cite{jwp}) have sometimes found it preferable 
to use a different and more sophisticated
approach, involving so--called ``dolittle diagrams'', to deal with
duality for Banach couples. A further discussion of these diagrams
is given on pp.\ 268--282 of the book \cite{bk} of Yuri Brudnyi
and Natan Krugljak.

\smallskip{}

We will slowly work our way through many well known facts. In particular
I will give a rather detailed proof of Calderón's famous duality theorem
for his complex interpolation spaces, using an approach with some
differences from his original proof. I believe that in some ways,
despite its considerable length, my proof of that theorem is simpler
than the original one, and draws attention to some hopefully useful
details about Calderón's interpolation spaces, which have not been
presented quite so explicitly in previous papers. Among the new features
in this second version of these notes are: 

(i) the mention (in Subsection \ref{sub:UsefulNormingSubspace}) of
an auxiliary result to Calderón's duality theorem which identifies
a useful norming subspace of the dual space of the space $\left[X_{0},X_{1}\right]_{\theta}$,
and

(ii) attention (in Definition \ref{def:3DifferentGs} and the first
part of the proof in Subsection \ref{sub:OurProof}) to the slight
difference between the spaces defined by Alberto Calder\'on and by
James Stafney which both of them denoted by $\mathcal{G}\left(X_{0},X_{1}\right)$
in their respective papers, and 

(iii) correction of a few small misprints in the first version (and
perhaps the introduction of some new ones!).

\smallskip{}

We will not make any explicit use here of the remarkable alternative
characterizations, due to Svante Janson \cite{janson}, of the complex
method of interpolation via weighted sequence spaces of Fourier coefficients.
These characterizations can be very useful for various purposes. They
provide an alternative way of describing the dual spaces of Calderón's
spaces. (Cf.\
\cite{janson} Theorems 1 and 2 on pp.\ 53--55 and, in particular,
Theorem 23 on p.\ 69.)

\smallskip{}

I apologize to those readers who might find many of the details mentioned
here to be superflous. But one of my aims is to make these notes comfortably
accessible to beginning graduate students, including those who might
perhaps be looking at my paper \cite{clat} about complex interpolation
of compact operators mapping into a couple of Banach lattices. At
the same time, ``conversely'', I should perhaps already make it
clear to the above-mentioned students, that these notes are not entirely
``self-contained''. In quite a number of places they will still
need to consult other books and papers for some essential details.

\smallskip{}

\section{\label{nbd}Notation, basic definitions and properties of the dual
couple of a Banach couple}

\smallskip{}
 All Banach spaces considered here will be assumed to be complex Banach
spaces.

For any two Banach spaces, $A$ and $B$ it will be convenient to
use the notation $A\overset{1}{\subset}B$ to mean that $A$ is contained
in $B$ and that, furthermore, $\left\Vert a\right\Vert _{B}\le\left\Vert a\right\Vert _{A}$
for each $a\in A$. The notation $A\overset{1}{=}B$ will mean that
$A\overset{1}{\subset}B$ and also $B\overset{1}{\subset}A$, i.e.,
it means that $A$ and $B$ coincide with equality of norms.

As usual, $A^{*}$ will denote the space of all continuous linear
functionals on $A$. But we will sometimes use the alternative notation
$A^{\#}$ for another space which is isometrically isomorphic to $A^{*}$.
\begin{defn}
\label{def:BanachPair}A \textbf{\textit{Banach couple}}, sometimes
also referred to as a \textbf{\textit{Banach pair}} or \textbf{\textit{interpolation
pair,}} and denoted by $(X_{0},X_{1})$, is an ordered pair of Banach
spaces $X_{0}$ and $X_{1}$ which are both contained in some Hausdorff
topological vector space $\mathcal{X}$. The two inclusion maps $X_{0}\subset\mathcal{X}$
and $X_{1}\subset\mathcal{X}$ are both required to be continuous. \end{defn}
\begin{rem}
In fact the conditions imposed in Definition \ref{def:BanachPair}
can be equivalently replaced by certain others which seem \textit{a
priori} to be weaker or stronger: 

First, with regard to weaker conditions, we shall see that we do not
really have to care about any topological structure of the vector
space $\mathcal{X}$ containing $X_{0}$ and $X_{1}$. It is enough
to require that the norms of $X_{0}$ and $X_{1}$ satisfy the so-called
``separability axiom'' which appears as (3.1) on p.~225 of \cite{PeetreJSparrG1972}
in the more general setting of normed abelian groups. I.e., it suffices
to require that whenever $y_{0}$ and $y_{1}$ are elements in $\mathcal{X}$
and whenever there exists a sequence $\left\{ x_{n}\right\} _{n\in\mathbb{N}}$
of elements in $\mathcal{X}$ such that, both for $j=0$ and $j=1$,
we have $x_{n}-y_{j}\in X_{j}$ for all $n$ and $\lim_{n\to}\left\Vert x_{n}-y_{j}\right\Vert _{X_{j}}=0$,
then $y_{0}=y_{1}$. In fact it is only necessary to require this
condition in the setting when $\mathcal{X}$ is replaced by its possibly
smaller subspace $X_{0}+X_{1}$. This suffices to show that the natural
seminorm on $X_{0}+X_{1}$ is a norm, such that $X_{0}$ and $X_{1}$
are both continuously embedded into $X_{0}+X_{1}$ equipped with this
norm. (Cf.~the latter part of the proof of Proposition 3.1 on p.~224
of \cite{PeetreJSparrG1972}.) 

As for imposing stronger conditions, we can also see, with the help
of the discussion in the preceding paragraph, that requiring the space
$\mathcal{X}$ in Definition \ref{def:BanachPair} to be normed, (or
even, see below, to be a Banach space) would give yet another equivalent
version of that definition.
\end{rem}
\begin{comment}
Suppose that we know that this separability condition holds. Then
the linear space $X{}_{0}+X_{1}$ exists as a linear subspace of $\mathcal{X}$.
Now define a seminorm $N(\cdot)$ on $X_{0}+X_{1}$ by setting 
\[
N(x):=\inf\left\{ \left\Vert x_{0}\right\Vert _{X_{0}}+\left\Vert x_{1}\right\Vert _{X_{1}}:x=x_{0}+x_{1},\mbox{ }x_{0}\in X_{0},\, x_{1}\in X_{1}\right\} .
\]
If we can show that this is a norm, then $X_{0}+X_{1}$ will be a
normed and therefore Hausdorff topological vectore space, which can
play the role of $\mathcal{X}$ in the above formal definition. Suppose
that $N(x)=0$. Then $x=u_{n}+v_{n}$ where $u_{n}\in X_{0}$ and
$v_{n}\in X_{1}$ and $\left\Vert u_{n}\right\Vert _{X_{0}}$ and
$\left\Vert v_{n}\right\Vert _{X_{1}}$ both tend to $0$ as $n\to\infty$.
Set $x_{n}=x-u_{n}=v_{n}$. Then $\lim_{n\to\infty}\left\Vert x_{n}-0\right\Vert _{X_{1}}=0$
and $\lim_{n\to\infty}\left\Vert x_{n}-x\right\Vert _{X_{0}}=0$ and
so the separability axiom gives us that $x=0$ completing the proof
that $N(\cdot)$ is a norm.
\end{comment}

For each Banach couple $(X_{0},X_{1})$, it is well known and immediately
evident that the spaces $X_{0}\cap X_{1}$ and $X_{0}+X_{1}$ normed,
respectively, by $\left\Vert x\right\Vert _{X_{0}\cap X_{1}}=\max\left\{ \left\Vert x\right\Vert _{X_{0}},\left\Vert x\right\Vert _{X_{1}}\right\} $
and $\left\Vert x\right\Vert _{X_{0}+X_{1}}=\inf\left\{ \left\Vert x_{0}\right\Vert _{X_{0}}+\left\Vert x_{1}\right\Vert _{X_{1}}:x=x_{0}+x_{1}\right\} $
are also Banach spaces.
\begin{defn}
\label{regi}Let $(X_{0},X_{1})$ be a Banach couple. We say that
$(X_{0},X_{1})$ is \textbf{\textit{regular}} if $X_{0}\cap X_{1}$
is dense in $X_{j}$ for $j=0,1$. A Banach space $X$ which satisfies
$X_{0}\cap X_{1}\subset X\subset X_{0}+X_{1}$ where both the inclusions
are continuous, is said to be an \textbf{\textit{intermediate space}}
with respect to the couple $(X_{0},X_{1})$. If, in addition, $X_{0}\cap X_{1}$
is dense in $X$ then we say that $X$ is a \textbf{\textit{regular
intermediate space}} with respect to $(X_{0},X_{1})$. 
\end{defn}
We shall assume throughout this section that $(X_{0},X_{1})$ is a
regular Banach couple of complex Banach spaces.

\smallskip{}
 The notation $\left\langle \cdot,\cdot\right\rangle $ will always
denote the bilinear functional defined on $(X_{0}\cap X_{1})\times(X_{0}\cap X_{1})^{*}$.
Thus, whenever we write $\left\langle x,y\right\rangle $, we will
always be assuming that $x\in X_{0}\cap X_{1}$ and $y\in(X_{0}\cap X_{1})^{*}$
and $\left\langle x,y\right\rangle $ is the value of the linear functional
$y$ acting on the element $x$. Later we shall define an explicit
variant of this notation for the case where $x$ is in a space larger
than $X_{0}\cap X_{1}$ and $y$ is in a suitable space smaller than
$(X_{0}\cap X_{1})^{*}$. 

It should be mentioned that, in later sections of these notes, it
will occasionally be convenient to use the notation $y(x)$ rather
than $\left\langle x,y\right\rangle $ or its variants, to denote
the value of a linear functional $y$ in the dual of some Banach space
when it acts on the element $x$ of that space. 
\begin{defn}
\label{blg}Suppose that $X$ is a regular intermediate space with
respect to the couple $(X_{0},X_{1})$. Let $X^{\#}$ denote the subspace
of elements $y\in(X_{0}\cap X_{1})^{*}$ which satisfy 
\[
\left\Vert y\right\Vert _{X^{\#}}:=\sup\left\{ \left|\left\langle x,y\right\rangle \right|:x\in X_{0}\cap X_{1},\left\Vert x\right\Vert _{X}\le1\right\} <\infty.
\]
 
\end{defn}
Note in particular that $(X_{0}\cap X_{1})^{\#}\overset{1}{=}(X_{0}\cap X_{1})^{*}$.
\begin{fact}
\label{fkz}Suppose that $U$ and $V$ are both regular intermediate
spaces with respect to $(X_{0},X_{1})$ and $U\overset{1}{\subset}V$.
Then $V^{\#}\overset{1}{\subset}U^{\#}$. More generally, if $U$
is continuously embedded in $V$ then $V^{\#}$ is continuously embedded
in $U^{\#}$. 
\end{fact}
\textit{Proof of Fact \ref{fkz}.} This will follow immediately from
the definition. We have $\left\Vert x\right\Vert _{V}\le c\left\Vert x\right\Vert _{U}$
for all $x\in X_{0}\cap X_{1}$ and some constant $c$ (and in particular
we will also consider the case where $c=1$). Then $\left|\left\langle x,y\right\rangle \right|\le\left\Vert x\right\Vert _{V}\left\Vert y\right\Vert _{V^{\#}}\le c\left\Vert x\right\Vert _{U}\left\Vert y\right\Vert _{V^{\#}}$
for all $x\in X_{0}\cap X_{1}$ and all $y\in V^{\#}$. Consequently
$y\in U^{\#}$ with $\left\Vert y\right\Vert _{U^{\#}}\le c\left\Vert y\right\Vert _{V^{\#}}$.
$\qed$

\smallskip{}
 In particular, Fact \ref{fkz} tells us that, for each intermediate
space $X$ in which $X_{0}\cap X_{1}$ is dense, we have that 
\begin{equation}
X^{\#}\subset(X_{0}\cap X_{1})^{\#}\overset{1}{=}(X_{0}\cap X_{1})^{*}\text{ and this inclusion is continuous.}\label{zaq}
\end{equation}

For $j=0,1$ we shall denote $X_{j}^{\#}=(X_{j})^{\#}$. Thus $X_{0}^{\#}\cap X_{1}^{\#}$
is the set of all elements $y\in(X_{0}\cap X_{1})^{*}$ which satisfy
\[
\left\Vert y\right\Vert _{X_{0}^{\#}\cap X_{1}^{\#}}=\max\left\{ \left\Vert y\right\Vert _{X_{0}^{\#}},\left\Vert y\right\Vert _{X_{1}^{\#}}\right\} <\infty.
\]

\smallskip{}
 Since $X_{0}\cap X_{1}$ is dense in $X_{0}+X_{1}$, the space $(X_{0}+X_{1})^{\#}$
can also be defined as in Definition \ref{blg} and we claim that
\begin{fact}
\label{kaler}$X_{0}^{\#}\cap X_{1}^{\#}$ coincides isometrically
with $(X_{0}+X_{1})^{\#}$. 
\end{fact}
\textit{Proof.} First, for each $x\in X_{0}\cap X_{1}$ and each $y\in X_{0}^{\#}\cap X_{1}^{\#}$
and each decomposition of $x$ in the form $x=x_{0}+x_{1}$, where
$x_{j}\in X_{j}$ (and so in fact $x_{j}\in X_{0}\cap X_{1}$) for
$j=0,1$, we have 
\begin{eqnarray*}
\left|\left\langle x,y\right\rangle \right| & = & \left|\left\langle x_{0},y\right\rangle +\left\langle x_{1},y\right\rangle \right|\le\left|\left\langle x_{0},y\right\rangle \right|+\left|\left\langle x_{1},y\right\rangle \right|\\
 & \le & \left\Vert x_{0}\right\Vert _{X_{0}}\left\Vert y\right\Vert _{X_{0}^{\#}}+\left\Vert x_{1}\right\Vert _{X_{1}}\left\Vert y\right\Vert _{X_{1}^{\#}}\\
 & \le & \left(\left\Vert x_{0}\right\Vert _{X_{0}}+\left\Vert x_{1}\right\Vert _{X_{1}}\right)\max\left\{ \left\Vert y\right\Vert _{X_{0}^{\#}},\left\Vert y\right\Vert _{X_{1}^{\#}}\right\} \\
 & = & \left(\left\Vert x_{0}\right\Vert _{X_{0}}+\left\Vert x_{1}\right\Vert _{X_{1}}\right)\left\Vert y\right\Vert _{X_{0}^{\#}\cap X_{1}^{\#}}.
\end{eqnarray*}
 Taking the infimum over all decompositions $x=x_{0}+x_{1}$ of the
kind specified above, we obtain that 
\[
\left|\left\langle x,y\right\rangle \right|\le\left\Vert x\right\Vert _{X_{0}+X_{1}}\left\Vert y\right\Vert _{X_{0}^{\#}\cap X_{1}^{\#}}.
\]
 Then, keeping $y$ fixed and taking the supremum over all $x\in X_{0}\cap X_{1}$
with $\left\Vert x\right\Vert _{X_{0}+X_{1}}\le1$, we obtain that
$X_{0}^{\#}\cap X_{1}^{\#}\overset{1}{\subset}(X_{0}+X_{1})^{\#}$.

Conversely, if $y\in(X_{0}+X_{1})^{\#}$, then, for $j=0,1$ we have,
for all $x\in X_{0}\cap X_{1}$, that 
\[
\left|\left\langle x,y\right\rangle \right|\le\left\Vert x\right\Vert _{X_{0}+X_{1}}\left\Vert y\right\Vert _{(X_{0}+X_{1})^{\#}}\le\left\Vert x\right\Vert _{X_{j}}\left\Vert y\right\Vert _{(X_{0}+X_{1})^{\#}}.
\]
 Taking the supremum over all $x\in X_{0}\cap X_{1}$ with $\left\Vert x\right\Vert _{X_{j}}\le1$,
we obtain that $y\in X_{j}^{\#}$ with $\left\Vert y\right\Vert _{X_{j}^{\#}}\le\left\Vert y\right\Vert _{(X_{0}+X_{1})^{\#}}$.
It follows that $(X_{0}+X_{1})^{\#}\overset{1}{\subset}X_{0}^{\#}\cap X_{1}^{\#}$.
So indeed we have 
\begin{equation}
(X_{0}+X_{1})^{\#}\overset{1}{=}X_{0}^{\#}\cap X_{1}^{\#}\label{ghg}
\end{equation}
 and the proof of Fact \ref{kaler} is complete. $\qed$

It seems intuitively obvious that, for each regular intermediate space
$X$, the space $X^{\#}$ can be identified with the dual space $X^{*}$.
But let us now explain this more precisely:
\begin{defn}
\label{ranX} Suppose that $X$ is a regular intermediate space with
respect to the regular Banach couple $(X_{0},X_{1})$.

(i) For each $y\in X^{\#}$ and each $x\in X$ we define 
\[
\left\langle x,y\right\rangle _{X}:=\lim_{n\rightarrow\infty}\left\langle x_{n},y\right\rangle 
\]
 where $\left\{ x_{n}\right\} _{n\in\mathbb{N}}$ is a sequence (any
sequence) of elements in $X_{0}\cap X_{1}$ which converges to $x$
in $X$ norm.

(ii) In contexts where it can help avoid ambiguity, we can use the
alternative and consistent notation $\left\langle \cdot,\cdot\right\rangle _{X{}_{0}\cap X_{1}}$
for the bilinear functional $\left\langle \cdot,\cdot\right\rangle $
introduced at the beginning of this discussion. =
\end{defn}
\smallskip{}

In many papers, notation such as $\left\langle x,y\right\rangle $,
originally introduced for the case when $x\in X_{0}\cap X_{1}$, is
still employed instead of $\left\langle x,y\right\rangle _{X}$, even
when $x$ is not in $X_{0}\cap X_{1}$, and apparently this usually
does not cause any problems or ambiguities. (This is done, for example,
on p.~1006 of \cite{ciddar}.) However here we feel the need to be
more explicit about the exact kind of duality being used at each step
of our development. (In Subsection \ref{sub:AnotherLook} below we
rewrite some of the arguments of the above-mentioned page 1006 of
\cite{ciddar} using the more pedantic terminology introduced here.) 

Of course the value of $\left\langle x,y\right\rangle _{X}$ is independent
of the particular choice of the sequence $\left\{ x_{n}\right\} _{n\in\mathbb{N}}$
and of course we also have 
\begin{equation}
\left|\left\langle x,y\right\rangle _{X}\right|\le\left\Vert x\right\Vert _{X}\left\Vert y\right\Vert _{X^{\#}}\text{ for each }x\in X\text{ and each }y\in X^{\#}.\label{qqw}
\end{equation}
It is also evident that $\left\langle \cdot,\cdot\right\rangle _{X}$
is a bilinear functional on $X\times X^{\#}$. Furthermore,
\begin{eqnarray}
\left\Vert y\right\Vert _{X^{\#}} & = & \sup\left\{ \left|\left\langle x,y\right\rangle \right|:x\in X_{0}\cap X_{1},\left\Vert x\right\Vert _{X}\le1\right\} \notag\\
 & = & \sup\left\{ \left|\left\langle x,y\right\rangle _{X}\right|:x\in X_{0}\cap X_{1},\left\Vert x\right\Vert _{X}\le1\right\} \notag\\
 & \le & \sup\left\{ \left|\left\langle x,y\right\rangle _{X}\right|:x\in X,\left\Vert x\right\Vert _{X}\le1\right\} \le\left\Vert y\right\Vert _{X^{\#}}\text{ ,}\label{aqqa}
\end{eqnarray}
 and so, in fact, all the inequalities in (\ref{aqqa}) are equalities. 

Let us use the notation $\lambda(x)$ for the action of a bounded
linear functional $\lambda\in X^{*}$ on an element $x\in X$. We
define the map $I_{X}:X^{\#}\rightarrow X^{*}$ by 
\[
(I_{X}y)(x)=\left\langle x,y\right\rangle _{X}\text{ for each }y\in X^{\#}\text{ and each }x\in X.
\]

It is clear from (\ref{aqqa}) that $\left\Vert I_{X}y\right\Vert _{X^{*}}=\left\Vert y\right\Vert _{X^{\#}}$
for each $y\in X^{\#}$.

\smallskip{}
 Now let $\lambda$ be an arbitrary element of $X^{*}$. Since $X$
is an intermediate space, for each $x\in X_{0}\cap X_{1}$ we have
$\left|\lambda(x)\right|\le\left\Vert \lambda\right\Vert _{X^{*}}\left\Vert x\right\Vert _{X}\le c\left\Vert \lambda\right\Vert _{X^{*}}\left\Vert x\right\Vert _{X_{0}\cap X_{1}}$
for some constant $c$. In other words, the restriction $\left.\lambda\right|_{X_{0}\cap X_{1}}$
of $\lambda$ to $X_{0}\cap X_{1}$ is a bounded linear functional
on $X_{0}\cap X_{1}$ and thus an element of $(X_{0}\cap X_{1})^{*}=(X_{0}\cap X_{1})^{\#}$.
We can thus use the notation 
\begin{equation}
\left\langle x,\left.\lambda\right|_{X_{0}\cap X_{1}}\right\rangle =\lambda(x)\text{ for all }x\in X_{0}\cap X_{1}.\label{ppf}
\end{equation}
 Furthermore, since $X_{0}\cap X_{1}$ is dense in $X$, the functionals
$\left.\lambda\right|_{X_{0}\cap X_{1}}$ and $\lambda$ are in one
to one correspondence, i.e., for two elements $\lambda$ and $\sigma$
of $X^{*}$, we have $\left.\lambda\right|_{X_{0}\cap X_{1}}=\left.\sigma\right|_{X_{0}\cap X_{1}}$
if and only if $\lambda=\sigma$. We deduce that the map $J_{X}:X^{*}\rightarrow(X_{0}\cap X_{1})^{*}$
defined by $J_{X}\lambda=\left.\lambda\right|_{X_{0}\cap X_{1}}$
and (\ref{ppf}) for each $\lambda\in X^{*}$ is one to one. Furthermore,
\[
\left|\left\langle x,J_{X}\lambda\right\rangle \right|=\left|\lambda(x)\right|\le\left\Vert \lambda\right\Vert _{X^{*}}\left\Vert x\right\Vert _{X}\text{ for all }x\in X_{0}\cap X_{1},
\]
 and so, for each $\lambda\in X^{*}$, we have $J_{X}\lambda\in X^{\#}$
with $\left\Vert J_{X}\lambda\right\Vert _{X^{\#}}\le\left\Vert \lambda\right\Vert _{X^{*}}$,
i.e., $J_{X}:X^{*}\rightarrow X^{\#}$ with norm not exceeding $1$.
It is clear that the operators $J_{X}$ and $I_{X}$ set up a one
to one isometry between $X^{*}$ and $X^{\#}$. Thus we have shown
that
\begin{fact}
\label{ning} $X^{\#}$ is exactly the dual space of $X$, if we adopt
the convention that the bilinear functional $\left\langle x,y\right\rangle _{X}$
defined in Definition \ref{ranX} specifies the action of each element
$y$ of the dual space on each element $x$ of $X$. 
\end{fact}
Now we can consider the couple consisting of the dual spaces $X_{0}^{\#}$
and $X_{1}^{\#}$.

\smallskip{}

\begin{fact}
\label{sumq}Let $(X_{0},X_{1})$ be a regular Banach couple. Then
$(X_{0}^{\#},X_{1}^{\#})$ is a Banach couple, and $(X_{0}\cap X_{1})^{\#}\overset{1}{=}X_{0}^{\#}+X_{1}^{\#}$. 
\end{fact}
\textit{Proof.} First we observe (cf.\ (\ref{zaq})) that $X_{0}^{\#}$
and $X_{1}^{\#}$ are both continuously embedded in $(X_{0}\cap X_{1})^{*}$.
So $(X_{0}^{\#},X_{1}^{\#})$ is indeed a Banach couple. More precisely,
since $X_{0}\cap X_{1}\overset{1}{\subset}X_{j}$, we have $X_{j}^{\#}\overset{1}{\subset}(X_{0}\cap X_{1})^{*}$
for $j=0,1$ and so $X_{0}^{\#}+X_{1}^{\#}\overset{1}{\subset}(X_{0}\cap X_{1})^{*}$.
As already mentioned, we have $(X_{0}\cap X_{1})^{\#}\overset{1}{=}(X_{0}\cap X_{1})^{*}$.
It remains to show that 
\begin{equation}
(X_{0}\cap X_{1})^{\#}\overset{1}{\subset}X_{0}^{\#}+X_{1}^{\#}.\label{aaau}
\end{equation}

Consider the Banach space $X_{0}\oplus_{\infty}X_{1}$\ which is
the cartesian product of $X_{0}$\ and $X_{1}$\ normed by $\left\Vert x_{0}\oplus x_{1}\right\Vert _{X_{0}\oplus_{\infty}X_{1}}:=\max\left\{ \left\Vert x_{0}\right\Vert _{X_{0}},\left\Vert x_{1}\right\Vert _{X_{1}}\right\} $.
Adopting the conventions which were fixed above in Definition \ref{ranX}
and the subsequent discussion for realizing the dual space of each
regular intermediate space with respect to the couple $(X_{0},X_{1})$,
we see that dual of $X_{0}\oplus_{\infty}X_{1}$\ is of course the
space $X_{0}^{\#}\oplus_{1}X_{1}^{\#}$, i.e., the cartesian product
of $X_{0}^{\#}$\ and $X_{1}^{\#}$\ with norm $\left\Vert y_{0}\oplus y_{1}\right\Vert _{X_{0}^{\#}\oplus_{1}X_{1}^{\#}}=\left\Vert y_{0}\right\Vert _{X_{0}^{\#}}+\left\Vert y_{1}\right\Vert _{X_{1}^{\#}}$,
and the action of linear functionals is given by the formula 
\[
(y_{0}\oplus y_{1})(x_{0}\oplus x_{1})=\left\langle x_{0},y_{0}\right\rangle _{X_{0}}+\left\langle x_{1},y_{1}\right\rangle _{X_{1}}.
\]

\smallskip{}
 Let $\Delta$ be the subspace of $X_{0}\oplus_{\infty}X_{1}$\ consisting of all elements of the
form $x\oplus x$\ for $x\in X_{0}\cap X_{1}$\ and equipped with
the norm of $X_{0}\oplus_{\infty}X_{1}$.

\smallskip{}
 Let $y$ be an arbitrary element of $(X_{0}\cap X_{1})^{\#}$. Then,
for each $x\in X_{0}\cap X_{1}$, we have $\left|\left\langle x,y\right\rangle \right|\le\left\Vert y\right\Vert _{(X_{0}\cap X_{1})^{\#}}\left\Vert x\right\Vert _{X_{0}\cap X_{1}}=\left\Vert y\right\Vert _{(X_{0}\cap X_{1})^{\#}}\max\left\{ \left\Vert x\right\Vert _{X_{0}},\left\Vert x\right\Vert _{X_{1}}\right\} $.
Thus the functional $\lambda$\ defined by $\lambda(x\oplus x)=\left\langle x,y\right\rangle $\ for
all $x\in X_{0}\cap X_{1}$\ defines a bounded linear functional
on $\Delta$ with norm not exceeding $\left\Vert y\right\Vert _{(X_{0}\cap X_{1})^{\#}}$.
By the Hahn--Banach theorem, there exists an element of $\left(X_{0}\oplus_{\infty}X_{1}\right)^{*}$\ which
is an extension of $\lambda$\ and has the same norm. I.e., there
exist $y_{0}\in X_{0}^{\#}$\ and $y_{1}\in X_{1}^{\#}$ such that
\begin{equation}
\left\langle x,y\right\rangle =\left\langle x,y_{0}\right\rangle _{X_{0}}+\left\langle x,y_{1}\right\rangle _{X_{1}}\label{rza}
\end{equation}
 for all $x\in X_{0}\cap X_{1}$\ and $\left\Vert y_{0}\right\Vert _{X_{0}^{\#}}+\left\Vert y_{1}\right\Vert _{X_{1}^{\#}}=\left\Vert \lambda\right\Vert _{\Delta^{*}}\le\left\Vert y\right\Vert _{(X_{0}\cap X_{1})^{\#}}$.
Since in (\ref{rza}) we are only considering elements $x$ in $X_{0}\cap X_{1}$,
we can rewrite (\ref{rza}) as $\left\langle x,y\right\rangle =\left\langle x,y_{0}\right\rangle +\left\langle x,y_{1}\right\rangle =\left\langle x,y_{0}+y_{1}\right\rangle $.
This shows that the element $y$ can be written in the form $y=y_{0}+y_{1}$,
where $y_{j}\in X_{j}^{\#}$ for $j=0,1$ and so $y\in X_{0}^{\#}+X_{1}^{\#}$
with 
\[
\left\Vert y\right\Vert _{X_{0}^{\#}+X_{1}^{\#}}\le\left\Vert y_{0}\right\Vert _{X_{0}^{\#}}+\left\Vert y_{1}\right\Vert _{X_{1}^{\#}}\le\left\Vert y\right\Vert _{(X_{0}\cap X_{1})^{\#}}.
\]
 This establishes (\ref{aaau}) and completes the proof of Fact \ref{sumq}.
\qed

\smallskip{}

\begin{fact}
\label{kjh}Let $X$ be a regular intermediate space with respect
to the regular Banach couple $(X_{0},X_{1})$. Then $X^{\#}$ is an
intermediate space with respect to the Banach couple $(X_{0}^{\#},X_{1}^{\#})$. 
\end{fact}
\textit{Proof.} In view of Fact \ref{fkz}, the continuous inclusions
$X_{0}\cap X_{1}\subset X\subset X_{0}+X_{1}$ imply the continuous
inclusions $(X_{0}+X_{1})^{\#}\subset X^{\#}\subset(X_{0}\cap X_{1})^{\#}$.
We complete the proof by applying Fact \ref{kaler} and Fact \ref{sumq}.
\qed

\smallskip{}

\begin{rem}
\label{intrs}Suppose that $U$ and $V$ are both regular intermediate
spaces with respect to the regular Banach couple $(X_{0},X_{1})$.
Then obviously both $U+V$ and $U\cap V$ are also intermediate spaces
with respect to $(X_{0},X_{1})$ and, clearly, $U+V$ is regular.
In all examples that we know, the space $U\cap V$ is also regular.
But it would be surprising if there are no examples where it is not
regular, and we invite the reader to produce such an example. In view
of our efforts long ago on the first two pages of \cite{ciddar} we
are particularly curious to know if such an example can be found,
or shown never to exist, in the special case where both $U$ and $V$
are complex interpolation spaces $U=[X_{0},X_{1}]_{\theta_{0}}$ and
$V=[X_{0},X_{1}]_{\theta_{1}}$ (of course with $\theta_{0}\ne\theta_{1}$). 
\end{rem}

\section{\label{cald}Calderón's theorem about duals of complex interpolation
spaces}

\subsection{\protect\label{sub:Preface} Preface}

In this section we will discuss complex interpolation spaces $[X_{0},X_{1}]_{\theta}$
as defined and studied by Alberto Calderón in \cite{ca}. For this
purpose we will also need to be familiar with a variant of these spaces,
denoted by $[X_{0},X_{1}]^{\theta}$ and obtained by a construction
which is sometimes referred to as Calderón's ``second'' or ``upper''
method. We refer to \cite{ca} Sections 2, 3, 5, 6 and 9.2 on pp.\
114--116 for the definitions of $[X_{0},X_{1}]_{\theta}$ and $[X_{0},X_{1}]^{\theta}$
and for the definitions%
\footnote{One should be careful to really formulate these definitions \textit{exactly}
as in \cite{ca}. As shown in \cite{cj}, a seemingly very minor change
in the definition of $\mathcal{F}\left(X_{0},X_{1}\right)$, as was
made perhaps unintentionally in some papers, can cause the space $\left[X_{0},X_{1}\right]_{\theta}$
to sometimes be different from that obtained when using the original
definition of $\mathcal{F}\left(X_{0},X_{1}\right)$.%
} of the spaces $\mathcal{F}(X_{0},X_{1})$ and $\mathcal{G}(X_{0},X_{1})$
and $\overline{\mathcal{F}}(X_{0},X_{1})$ of Banach space valued
analytic functions which are needed to define and study $[X_{0},X_{1}]_{\theta}$
and $[X_{0},X_{1}]^{\theta}$. 

Some confusion could ultimately result from the fact that the notation
$\mathcal{G}(X_{0},X_{1})$ refers to (at least) three different spaces,
depending on which book or paper you happen to consult regarding complex
interpolation. Therefore, for our purposes, we will try to avoid such
confusion by the following somewhat informally expressed

\begin{defn}
\label{def:3DifferentGs}Let $\mathcal{G}_{\mathbf{Ca}}(X_{0},X_{1})$,
$\mathcal{G}_{\mathbf{St}}(X_{0},X_{1})$ and $\mathcal{G}_{\mathbf{BL}}(X_{0},X_{1})$
denote the three different spaces which are each denoted by $\mathcal{G}(X_{0},X_{1})$
by Calder\'on on p.~116 of \cite{ca}, by Stafney on p.~334 of
\cite{stafney}, and by Bergh-L\"ofstr\"om on pp.~88--89 of \cite{bl}
, respectively.
\end{defn}
It will be convenient to use the notation $\mathbb{S}$ for the ``unit
strip'' in the complex plane, i.e., 
\[
\mathbb{S}=\left\{ z\in\mathbb{C}:0\le\mathrm{Re\,}z\le1\right\} 
\]
 and to denote the interior of $\mathbb{S}$ by $\mathbb{S}^{\circ}$.
The boundary $\partial\mathbb{S}$ of $\mathbb{S}$ is of course the
union of the two vertical lines $\{it:t\in\mathbb{R}\}$ and $\{1+it:t\in\mathbb{R}\}$.

\smallskip{}
 We recall (\cite{ca} Sections 9.2--9.3 p.\ 116 and Sections 29.2--29.3
pp.\ 132--134) that $X_{0}\cap X_{1}$ is dense in $[X_{0},X_{1}]_{\theta}$
for each $\theta\in(0,1)$. Thus the dual space of $[X_{0},X_{1}]_{\theta}$
coincides with the space $\left([X_{0},X_{1}]_{\theta}\right)^{\#}$
defined as in Section \ref{nbd}. Among other basic properties of
$\left[X_{0},X_{1}\right]_{\theta}$ that we will need, we mention
the easily established fact that 
\begin{equation}
X_{0}\cap X_{1}\overset{1}{\subset}\left[X_{0},X_{1}\right]_{\theta}\overset{1}{\subset}X_{0}+X_{1}.\label{eq:IntSp}
\end{equation}

In the context and the notation that we have set up in Section \ref{nbd},
in particular bearing Fact \ref{ning} in mind, Calderón's remarkable
duality theorem (\cite{ca} Section 12.1 p.\ 121) for his spaces
$[X_{0},X_{1}]_{\theta}$ can be expressed as follows:
\begin{thm}
\label{crt}Let $(X_{0},X_{1})$ be a regular Banach couple. Then,
for each $\theta\in(0,1)$, the dual space $([X_{0},X_{1}]_{\theta})^{\#}$
coincides isometrically with $[X_{0}^{\#},X_{1}^{\#}]^{\theta}$. 
\end{thm}
\smallskip{}
 In other words, the element $y\in(X_{0}\cap X_{1})^{\#}$ satisfies
\[
\left\Vert y\right\Vert _{([X_{0},X_{1}]_{\theta})^{\#}}:=\sup\left\{ \left|\left\langle x,y\right\rangle \right|:x\in X_{0}\cap X_{1},\left\Vert x\right\Vert _{[X_{0},X_{1}]_{\theta}}\le1\right\} <\infty
\]
 if and only if $y$ is also an element of $[X_{0}^{\#},X_{1}^{\#}]^{\theta}$.
Furthermore, in that case, 
\[
\left\Vert y\right\Vert _{([X_{0},X_{1}]_{\theta})^{\#}}=\left\Vert y\right\Vert _{[X_{0}^{\#},X_{1}^{\#}]^{\theta}}.
\]

\smallskip{}

\begin{rem}
\label{rem:IfNotRegular}Of course there is a simple procedure to
make this theorem immediately applicable also to every Banach couple
$\left(A_{0},A_{1}\right)$ of complex Banach spaces. If such a couple
is not regular, we can let $X_{j}$ be the closure of $A_{0}\cap A_{1}$
in $A_{j}$ for $j=0,1$, thus obtaining the regular couple $\left(X_{0},X_{1}\right)$.
We will refer to the couple $\left(X_{0},X_{1}\right)$ obtained in
this way as the \textit{regularization} of $\left(A_{0},A_{1}\right)$.
We shall use the fact (see \cite{ca} p.~116 Section 9.3 p.~116
and Section 29.3 pp.~133--134) that $\left[A_{0},A_{1}\right]_{\theta}\overset{1}{=}\left[X_{0},X_{1}\right]_{\theta}$. 
\end{rem}

\begin{rem}
\label{rem:WeakerAnalyticity}Some relevant definitions and basic
properties of analytic Banach space valued functions are recalled,
for anyone who may need them, in an appendix (Subsection\ref{hilphi}).
Since the space $\left[X_{0}^{\#},X_{1}^{\#}\right]^{\theta}$ is
defined via the space $\overline{\mathcal{F}}(X_{0}^{\#},X_{1}^{\#})$
we should make some remarks about $\overline{\mathcal{F}}(X_{0}^{\#},X_{1}^{\#})$.
In the definition (see Section 5 on p.~115 of \cite{ca}) of the
space $\overline{\mathcal{F}}(B_{0},B_{1})$ for an arbitrary Banach
couple $\left(B_{0},B_{1}\right)$, one of the requirements for a
function $f:\mathbb{S}\to B_{0}+B_{1}$ to be in $\overline{\mathcal{F}}(B_{0},B_{1})$
is that its restriction to $\mathbb{S}^{\circ}$ must be an analytic
$B_{0}+B_{1}$ valued function on $\mathbb{S}^{\circ}$. (Cf.~Definition
\ref{def:HolomVecValued}). In cases where, as in the setting of Theorem
\ref{crt}, $\left(B_{0},B_{1}\right)$ is the ``dual'' couple $\left(X_{0}^{\#},X_{1}^{\#}\right)$
of a regular couple $\left(X_{0},X_{1}\right)$, one convenient way
(among several other equivalent ways) for expressing this analyticity
requirement is the condition that the function $\zeta\mapsto\left\langle x,f(\zeta)\right\rangle $
must be a scalar holomorphic function on $\mathbb{S}^{\circ}$ for
each element $x$ in the predual $X_{0}\cap X_{1}$ of $X_{0}^{\#}+X_{1}^{\#}$.
As explained/recalled in Subsection \ref{hilphi}, this condition
suffices to imply rather stronger conditions on $f$, among them the
existence of the complex derivative $f^{\prime}:\mathbb{S}^{\circ}\to X_{0}^{\#}+X_{1}^{\#}$,
which necessarily satisfies $\lim_{h\to0}\left\Vert \frac{1}{h}\left(f\left(\zeta+h\right)-f(\zeta)\right)-f'(\zeta)\right\Vert _{X_{0}^{\#}+X_{1}^{\#}}=0$
and therefore also 
\begin{equation}
\left\langle x,f^{\prime}(\zeta)\right\rangle =\frac{d}{d\zeta}\left\langle x,f(\zeta)\right\rangle \mbox{ for each }\zeta\in\mathbb{S}^{\circ}\mbox{ and each }x\in X_{0}\cap X_{1}.\label{eq:vvdia}
\end{equation}

\smallskip{}

\end{rem}
Our main goal in this section will be to give an alternative proof
of Theorem \ref{crt}. In keeping with the ``pedantic'' spirit of
these notes, it will use the notions and notation of Section \ref{nbd}.
Of course we are inspired by Calderón's beautiful and ingenious proof
in \cite{ca} Section 32.1 pp.\ 148--156 and by the similar proof%
\footnote{As already mentioned above, some of the notation used in \cite{bl}
is slightly different from that in \cite{ca}. In particular, $\overline{\mathcal{F}}\left(X_{0},X_{1}\right)$
is denoted by $\mathcal{G}\left(X_{0},X_{1}\right)$.%
}presented in \cite{bl} pp.\ 98--101. But perhaps in some ways, despite
its considerable length, our proof here might be considered to be
simpler.

\smallskip{}
Here are some features of our proof:

$\bullet$ In \cite{ca} and in \cite{bl} the inclusion $[X_{0}^{\#},X_{1}^{\#}]^{\theta}\overset{1}{\subset}([X_{0},X_{1}]_{\theta})^{\#}$
is proved using a somewhat elaborate multilinear interpolation theorem
(Section 11.1 p.\ 120 and Section 31.1 on pp 142--148 of \cite{ca}
and Theorem 4.4.2 on p.\ 97 of \cite{bl}.) Here we will be able
to manage without this auxiliary theorem. However, instead, we will
have to use an apparently simpler (and very attractive!) lemma of
Stafney.

$\bullet$ Calderón uses a connection between analytic functions on
$\mathbb{S}^{\circ}$ and analytic functions on the open unit disk.
In our proof the transition to analytic functions on the open unit
disk will be made at an earlier stage. Some readers might like the
fact that this approach saves us from having to deal with the Poisson
kernels for $\mathbb{S}^{\circ}$ which are given by somewhat complicated
formul\ae\ (which can be seen, for example, on p.\ 117 of \cite{ca}).
Indeed we have set ourselves the goal of using only the simplest possible
forms of various results about scalar valued analytic functions. Thus,
for example, we will essentially avoid having to talk about boundary
behaviour of bounded analytic functions, nontangential limits, properties
of the boundary values of harmonic conjugates of harmonic functions
with smooth boundary values, etc.

$\bullet$ We will use an alternative and in some ways more convenient
definition of the space $[Y_{0},Y_{1}]^{\theta}$ which is available
in the particular case, the only case that we need here, where $(Y_{0},Y_{1})=(X_{0}^{\#},X_{1}^{\#})$,
i.e., where $(Y_{0},Y_{1})$ is a couple of dual spaces of some regular
couple $(X_{0},X_{1})$. This alternative definition is analogous
to the definition of the space denoted by $B^{+}\{z_{0}\}^{*}$ in
\cite{ccrsw} pp.\ 211--212 and indeed some parts of our proof here
are suggested by things in \cite{ccrsw}. (In the setting of \cite{ccrsw}
where one deals with $n$--tuples or infinite families of Banach spaces
instead of couples, there is apparently no alternative to using this
kind of alternative definition of an ``upper'' method.)

$\bullet$ One consequence of using the alternative definition of
$[Y_{0},Y_{1}]^{\theta}$ mentioned just above is that there will
be no need to use the special duality between spaces $L^{1}(A)$ and
$\Lambda^{\prime}(A^{*})$ of Banach space valued functions, which
is defined and discussed somewhat elaborately in pp.~148--151 of
\cite{ca} (and more concisely on p.~98 and p.~101 of \cite{bl})
via appropriately generalized Riemann-Stieltjes integrals $\int_{-\infty}^{\infty}\left\langle f(t),dg(t)\right\rangle $
and considerations of density. Instead of this, we will only need
the standard duality between the spaces $L^{1}$ and $L^{\infty}$
of complex valued functions on $\mathbb{T}$.
\begin{rem}
We refer to pp.~227--231 of \cite{KreinSPetuninJSemenovE1982} for
yet another approach to describing the dual of $\left[X_{0},X_{1}\right]_{\theta}$.
There the dual is shown to coincide isometrically with $\widetilde{\left[X_{0}^{*},X_{1}^{*}\right]_{\theta}}$,
the Gagliardo completion of $\left[X_{0}^{*},X_{1}^{*}\right]_{\theta}$.
In general, for any Banach couple $\left(A_{0},A_{1}\right)$ the
Gagliardo completion of $\left[A_{0},A_{1}\right]_{\theta}$ is the
Banach space $\widetilde{\left[A_{0},A_{1}\right]_{\theta}}$ whose
unit ball is the closure of the unit ball of $\left[A_{0},A_{1}\right]_{\theta}$
with respect to the norm of $A_{0}+A_{1}$. It is not difficult to
show that $\left[A_{0},A_{1}\right]^{\theta}\overset{1}{\subset}\widetilde{\left[A_{0},A_{1}\right]_{\theta}}$
(with the help of (\ref{eq:NormCgce2Deriv}) and of sequences in $\mathcal{F}\left(A_{0},A_{1}\right)$
like the one introduced at the beginning of the proof of Theorem \ref{babaeq}).
We do not know of an example of a couple $\left(A_{0},A_{1}\right)$
for which this inclusion is strict. By combining the result of \cite[p. 227-231]{KreinSPetuninJSemenovE1982}
with the Calder\'on duality theorem, we see that it is an isometric
equality, at least in all those cases where $\left(A_{0},A_{1}\right)$
is the dual couple $\left(X_{0}^{\#},X_{1}^{\#}\right)$ of some couple
$\left(X_{0},X_{1}\right)$.
\end{rem}

\subsection{\label{svaf}Preliminaries, mainly about scalar valued analytic functions}

\smallskip{}

A number of rather standard results about scalar valued harmonic and
analytic functions will be used in our proof. I prefer that we discuss
them in this subsection, in advance, and in considerable detail, before
we start that proof. There are also two simple results about Banach
space valued analytic functions which I will include here. 

This is one place where I might particularly recall my apology at
the beginning of these notes. 

As usual $\mathbb{T}$ denotes the unit circle, $\mathbb{T}=\left\{ w\in\mathbb{C}:\left|w\right|=1\right\} $.
We also let $\mathbb{D}$ denote the closed unit disk $\mathbb{D}=\left\{ w\in\mathbb{C}:\left|w\right|\le1\right\} $
and of course $\mathbb{D}^{\circ}$ denotes its interior.

\smallskip{}

\begin{comment}
\textit{\textcolor{red}{Maybe in future versions work with $\mathbb{S}\cup\left\{ \infty\right\} $
instead of $\mathbb{S}$.}}\textit{\textcolor{red}{\small Would this
simplify the presentation?}}
\end{comment}
There are of course conformal maps of $\mathbb{D}^{\circ}$ onto $\mathbb{S}^{\circ}$.
We will find it convenient to use a particular one of these, which
we will denote by $\xi:\mathbb{D}^{\circ}\rightarrow\mathbb{S}^{\circ}$,
and whose choice depends on the value of the parameter $\theta\in(0,1)$
appearing in the statement of Theorem \ref{crt}. In fact we will
define $\xi$ on a larger set than $\mathbb{D}^{\circ}$ and we will
also need to take note of the behaviour of $\xi$ and its inverse
mapping $\xi^{-1}$ on the boundaries of $\mathbb{D}$ and $\mathbb{S}$
respectively.

As the first step of the (rather standard) construction of $\xi$
we consider the Möbius transformation 
\begin{equation}
\mu(w)=\frac{we^{-i\pi\theta}-e^{i\pi\theta}}{w-1},`\label{defmu}
\end{equation}
which maps $1$ to $\infty$, $e^{i2\pi\theta}$ to $0$, and each
other point on the unit circle $\mathbb{T}=\{w\in\mathbb{C}:$ $\left|w\right|=1\}$
to a non zero real number. More explicitly, each $e^{it}\in\mathbb{T}\backslash\{1,e^{i2\pi\theta}\}$
is mapped to 
\begin{equation}
\mu(e^{it})=\frac{e^{i(t-\pi\theta)}-e^{i\pi\theta}}{e^{it}-1}=\frac{e^{i(t/2-\pi\theta)}-e^{i(\pi\theta-t/2)}}{e^{it/2}-e^{-it/2}}=\frac{\sin(t/2-\pi\theta)}{\sin t/2}.\label{mdia}
\end{equation}

\smallskip{}
We introduce the two subintervals $I_{0}:=(2\pi\theta,2\pi)$ and
$I_{1}:=(0,2\pi\theta)$ of $(0,2\pi)$, and the two arcs $\Gamma_{0}$
and $\Gamma_{1}$ defined by 
\begin{equation}
\Gamma_{j}:=\left\{ e^{it}:t\in I_{j}\right\} \mbox{ for }j=0\mbox{ and }j=1.\label{eq:defgamj}
\end{equation}
We see from (\ref{mdia}) that $\mu$ maps $\Gamma_{0}$ onto the
positive real axis, and $\Gamma_{1}$ onto the negative real axis.

We also have $\mu(0)=e^{i\pi\theta}$. Thus, since $\theta\in(0,1)$,
the transformation $\mu$ conformally maps the open unit disk $\mathbb{D}^{\circ}=\left\{ w\in\mathbb{C}:\left|w\right|<1\right\} $
onto the open upper half plane and is analytic and conformal at every
point $w\ne1$ of the entire complex plane. 

As in \cite{ca}, this leads us to introduce the map $\xi$ defined
by 
\begin{equation}
\xi(w)=\frac{1}{i\pi}\log\mu(w)=\frac{1}{i\pi}\log\left(\frac{we^{-i\pi\theta}-e^{i\pi\theta}}{w-1}\right),\label{defxi}
\end{equation}
where here $\log z$ denotes a branch of the complex logarithm which
is analytic and conformal at every point of the punctured closed half-plane
$\Omega=\left\{ z:\mathrm{Re}\, z\ge0\right\} \setminus\left\{ 0\right\} $
and which satisfies $\mathrm{0\le Im}\,\log z\le\pi$ for all $z\in\Omega$.
It follows that $\xi$ maps $\mathbb{D}^{\circ}$ conformally onto
$\mathbb{S}^{\circ}$ and, furthermore, is also analytic and conformal
at every point of $\mathbb{D}\backslash\{1,e^{i2\pi\theta}\}$. We
note that $\xi(0)=\theta$ and that the images of the arcs $\Gamma_{0}$
and $\Gamma_{1}$ satisfy 
\begin{equation}
\xi(\Gamma_{j})=\{j+is:s\in\mathbb{R}\}\text{ for }j=0,1.\label{idrvba}
\end{equation}
Thus $\xi$ maps $\mathbb{D}\backslash\{1,e^{i2\pi\theta}\}$ continuously
onto $\mathbb{S}$. If $w$ tends to $1$ or to $e^{i2\pi\theta}$
then $\left|\mu(w)\right|$ tends to $\infty$ or to $0$ and therefore
$\left|\mathrm{Im}\,\xi(w)\right|$ tends to $\infty$. We shall also
need to know that $\xi$ is a one-to-one map, not just on $\mathbb{D}^{\circ}$,
but also on the whole of $\mathbb{D}\backslash\{1,e^{i2\pi\theta}\}$.
This will be an obvious consequence of the existence of its inverse
map $\xi^{-1}$, which we shall now describe: 

To obtain that inverse map, we note that, for each $z\in\mathbb{S}$
and each $w\in\mathbb{D}\setminus\left\{ 1,e^{i2\pi\theta}\right\} $,
if $\xi(w)=z$, then $\frac{we^{-i\pi\theta}-e^{i\pi\theta}}{w-1}=e^{i\pi z}$
and, consequently, 
\begin{equation}
w=\xi^{-1}(z):=\frac{e^{i\pi\theta}-e^{i\pi z}}{e^{-i\pi\theta}-e^{i\pi z}}.\label{eq:defxiinv}
\end{equation}
So we see that the function $\xi^{-1}$ is the restriction to $\mathbb{S}$
of a function which is analytic at every point $z\in\mathbb{C}$ which
is not an element of the sequence $\{2k-\theta\}_{k\in\mathbb{Z}}$.
More relevantly for us, $\xi^{-1}$ maps $\mathbb{S}^{\circ}$ conformally
onto $\mathbb{D}^{\circ}$, and is also conformal at every point of
$\partial\mathbb{S}$. Thus it is also a continuous map of $\mathbb{S}$
onto $\mathbb{D}\setminus\left\{ 1,e^{i\pi\theta}\right\} $, and
also of course one-to-one. In view of (\ref{idrvba}), it maps $\left\{ j+it:t\in\mathbb{R}\right\} $
onto $\Gamma_{j}$ for $j=0,1$. Furthermore, we note that $\xi^{-1}(z)$
tends to $e^{i2\pi\theta}$ or to $1$ as the imaginary part of the
point $z\in\mathbb{S}$ tends to $+\infty$ or, respectively, to $-\infty$.
More precisely, writing $z=x+iy$, we have 
\begin{equation}
\lim_{R\to+\infty}\sup\left\{ \left|\xi^{-1}(x+iy)-e^{i2\pi\theta}\right|:x\in[0,1],\, y>R\right\} =0\label{eq:upstrip}
\end{equation}
and 
\begin{equation}
\lim_{R\to+\infty}\sup\left\{ \left|\xi^{-1}(x+iy)-1\right|:x\in[0,1],\, y<-R\right\} =0,\label{eq:DownStrip}
\end{equation}
these being consequences of the following simple estimates: 
\begin{eqnarray*}
\left|\frac{e^{i\pi\theta}-e^{i\pi z}}{e^{-i\pi\theta}-e^{i\pi z}}-e^{i2\pi\theta}\right| & = & \left|\frac{-e^{i\pi z}+e^{i\pi(2\theta+z)}}{e^{-i\pi\theta}-e^{i\pi z}}\right|\le\frac{\left|e^{i\pi z}\left(-1+e^{i2\pi\theta}\right)\right|}{\left|e^{-i\pi\theta}\left(1-e^{i\pi(z-\theta)}\right)\right|}\\
 & \le & \frac{2\left|e^{i\pi z}\right|}{\left|\left|e^{i\pi(z-\theta)}\right|-1\right|}=\frac{2e^{-\pi y}}{\left|\left|e^{-\pi y}\right|-1\right|}.
\end{eqnarray*}
and 
\begin{eqnarray*}
\left|\frac{e^{i\pi\theta}-e^{i\pi z}}{e^{-i\pi\theta}-e^{i\pi z}}-1\right| & = & \left|\frac{2i\sin\pi\theta}{e^{-i\pi\theta}-e^{i\pi z}}\right|\le\frac{2}{\left|e^{-i\pi\theta}\left(1-e^{i\pi(z-\theta)}\right)\right|}\\
 & \le & \frac{2}{\left|\left|e^{i\pi(z-\theta)}\right|-1\right|}=\frac{1}{\left|\left|e^{-\pi y}\right|-1\right|}.
\end{eqnarray*}

\smallskip{}

Later on, we will need to use the following very simple lemma, in
which the functions $\xi$ and $\xi^{-1}$ give us convenient ways
of making transitions between certain functions $\phi$ defined on
$\mathbb{S}$ and certain other functions $\psi$ defined on $\mathbb{D}$.
It is convenient for our needs, and just as easy, to obtain this lemma,
not just for scalar functions but also for analytic Banach space valued
functions. 
\begin{lem}
\label{lem:PhiPsi}Let $A$ be an arbitrary Banach space. 

(i) Suppose that $\phi:\mathbb{S}\to A$ be an arbitrary continuous
function which is analytic in $\mathbb{S}^{\circ}$ and which vanishes
at $\infty$, i.e., which satsfies 
\begin{equation}
\lim_{R\to+\infty}\sup\left\{ \left\Vert \phi(z)\right\Vert _{A}:z\in\mathbb{S},\,\left|\mathrm{Im}\, z\right|>R\right\} =0.\label{eq:vai}
\end{equation}
Then the function $\psi:\mathbb{D}\to A$ defined by $\psi(1)=\psi(e^{i2\pi\theta})=0$
and $\psi(z)=\phi(\xi(z))$ for all $z\in\mathbb{D}\setminus\{1,e^{i2\pi\theta}\}$
is continuous on $\mathbb{D}$ and analytic in $\mathbb{D}^{\circ}$.

(ii) Conversely, if $\psi:\mathbb{D}\to A$ is an arbitrary continuous
function which is analytic in $\mathbb{D}^{\circ}$ and satisfies
$\psi(1)=\psi(e^{i2\pi\theta})=0$, then the function $\phi:\mathbb{S}\to A$
defined by $\phi(z)=\psi(\xi^{-1}(z))$ is continuous on $\mathbb{S}$,
analytic in $\mathbb{S}^{\circ}$ and satisfies (\ref{eq:vai}). 
\end{lem}
\noindent \textit{Proof.} This is an obvious exercise, using the
various properties of $\xi$ and $\xi^{-1}$ which we have established
just above. The case of general $A$ is essentially the same as the
case $A=\mathbb{C}$. (Recall Definition \ref{def:HolomVecValued},
where of course one can take $Y=A^{*}$.) $\qed$ %
\begin{comment}
\textit{\textcolor{red}{Maybe be more explicit about this in future
versions. Maybe upgrade to a more elaborate version which will also
take care more explicitly of the claim after Definition \ref{deffd}.
Maybe put details in an appendix.}}
\end{comment}

\smallskip{}

We recall the following standard result of Phragmen--Lindelöf type,
which is often referred to as the ``three lines theorem'' and associated
with the names of Hadamard and Doetsch. (More precisely, the analogous
``three circles theorem'' is due to Jacques Hadamard and its variant
for three lines is due to Gustav Doetsch.) %% Doetsch apparently had a Nazi connection.
We state it here in a version which is adequate for our particular
needs. (In fact (\ref{jtd}) can be replaced by a very much weaker
condition.)
\begin{lem}
\label{3l}Let $M_{0}$, $M_{1}$ and $c$ be positive constants.

(i) Suppose that $f:\mathbb{S}\rightarrow\mathbb{C}$ is a continuous
function which is analytic in $\mathbb{S}^{\circ}$ and satisfies
\begin{equation}
\left|f(z)\right|\le c(1+\left|z\right|)\text{ for all }z\in\mathbb{S}^{\circ}\label{jtd}
\end{equation}
 and 
\begin{equation}
\sup_{t\in\mathbb{R}}\left|f(j+it)\right|\le M_{j}\text{ for }j=0,1.\label{wqz}
\end{equation}
 Then 
\begin{equation}
\left|f(z)\right|\le M_{0}^{1-\mathrm{Re\,}z}M_{1}^{\mathrm{Re\,}z}\text{ for all }z\in\mathbb{S}^{\circ}.\label{pqwz}
\end{equation}

(ii) The same conclusion (\ref{pqwz}) holds if $f$ is not defined
on $\partial\mathbb{S}$ and is merely an analytic function on $\mathbb{S}^{\circ}$
satisfying (\ref{jtd}) and which, instead of condition (\ref{wqz}),
satisfies 
\begin{equation}
\lim_{r\searrow0}\sup\left\{ \left|f(z)\right|:z\in\mathbb{S},\text{ }0<\left|\mathrm{Re\,}z-j\right|<r\right\} \le M_{j}\text{ for }j=0,1.\label{swqz}
\end{equation}
 
\end{lem}
\textit{Proof.} It is almost quicker to give a proof than a reference,
and, anyway, parts of the proof here will also be relevant later.

We first deal with part (i): Let $\delta$ be a positive constant
and consider the entire function $\phi(z)=\exp(\delta z^{2}-(1-z)\log M_{0}-z\log M_{1})$.
Fix $z_{*}\in\mathbb{S}^{\circ}$. By the maximum modulus principle,
$\left|\phi(z_{*})f(z_{*})\right|\le\sup_{z\in\partial E_{R}}\left|\phi(z)f(z)\right|$
where $E_{R}$ is the rectangle $\{z\in\mathbb{S}:\left|\mathrm{Im\,}z\right|\le R\}$
and $R>\left|z_{*}\right|$. In view of \smallskip{}

\begin{equation}
\left|e^{\delta z^{2}}\right|=e^{\delta[(\mathrm{Re\,}z)^{2}-(\mathrm{Im\,}z)^{2}]}\label{uglfe}
\end{equation}
 and (\ref{wqz}) we see that $\left|\phi(j+it)f(j+it)\right|\le e^{\delta j}$
for $j=0,1$ and all $t\in\mathbb{R}$ and, in view of (\ref{jtd})
and (\ref{uglfe}), we can also ensure, by choosing $R$ sufficiently
large, that $\left|\phi(z)f(z)\right|\le1$ at every point $z$ on
each of the two horizontal line segments of $\partial E_{R}$. We
deduce that $\left|\phi(z_{*})f(z_{*})\right|\le e^{\delta}$. So
$\left|f(z_{*})\right|\le e^{\delta}/|\phi(z_{*})|=e^{\delta(1-z_{*}^{2})}|M_{0}^{1-\mathrm{Re\,}z_{*}}M_{1}^{\mathrm{Re\,}z_{*}}$.
Since we can choose $\delta$ arbitrarily small, this implies (\ref{pqwz}).

\smallskip{}
 Part (ii) is a straightforward corollary of part (i). In this case,
given an analytic function $f:\mathbb{S}^{\circ}\rightarrow\mathbb{C}$
which satisfies (\ref{jtd}) and (\ref{swqz}), we simply apply part
(i) to the function $f_{r}(z):=f\left(r/2+(1-r)z\right)$ where $r\in(0,1/2)$
is a constant. It is clear that $f_{r}$ satisfies all the hypotheses
of part (i), with $c$ replaced by some other constant $c_{r}$ depending
on $r$, and $M_{j}$ replaced by $\sup_{t\in\mathbb{R}}\left|f(|j-r/2|+it)\right|$.
Thus part (i) gives us, for each $z\in\mathbb{S}^{\circ}$, that 
\[
\left|f_{r}(z)\right|\le\left(\sup\left\{ \left|f(\zeta)\right|:0<\mathrm{Re\,}\zeta<r\right\} \right)^{1-\mathrm{Re\,}z}\cdot\left(\sup\left\{ \left|f(\zeta)\right|:0<1-\mathrm{Re\,}\zeta<r\right\} \right)^{\mathrm{Re\,}z}
\]
 and we obtain (\ref{pqwz}) for $f$ by passing to the limit as $r$
tends to $0$. \qed
\begin{rem}
A standard example%
\footnote{For example, consider the function $f(z)=\exp\left(\exp\left(i\pi(z-1/2)\right)\right)$,
which is an obvious modification of the function discussed on page
256 of \cite{rudin}.%
} shows that the conclusions of Lemma \ref{3l} may fail to hold if
the condition (\ref{jtd}) is not imposed. However, as already mentioned,
(\ref{jtd}) can be replaced by much weaker conditions. 
\end{rem}
\smallskip{}

This next result could be considered to be a variant of Lemma \ref{3l},
in which $M_{0}=M_{1}$ and we only consider the behaviour of the
function ``near infinity''. It is convenient for our later needs
to allow that function to take values in a Banach space. The proof
is only a minor variant of the proof needed for the case of scalar
valued functions. In fact we can already state those ``later needs''
as an immediate corollary of the lemma.
\begin{lem}
\label{lem:3LinesAtInfty}Let $A$ be a Banach space, whose norm will
be denoted by $\left\Vert \cdot\right\Vert $. Suppose that the function
$f:\mathbb{S}\to A$ is continuous and bounded on $\mathbb{S}$, and
analytic in $\mathbb{S}^{\circ}$. Suppose further that, for $j=0$
and $j=1$, $\lim_{t\to-\infty}\left\Vert f(j+it)\right\Vert =\lim_{t\to+\infty}\left\Vert f(j+it)\right\Vert =0$.
Then 
\begin{equation}
\lim_{R\to\infty}\sup\left\{ \left\Vert f(z)\right\Vert :z\in\mathbb{S},\,\left|\mathrm{Im}\, z\right|\ge R\right\} =0.\label{eq:awe}
\end{equation}
\end{lem}
\begin{cor}
\label{cor:VanishAtInfty}For each Banach couple $\left(X_{0},X_{1}\right)$,
each function $f\in\mathcal{F}\left(X_{0},X_{1}\right)$ has the property
\[
\lim_{R\to\infty}\sup\left\{ \left\Vert f(z)\right\Vert _{X_{0}+X_{1}}:z\in\mathbb{S},\,\left|\mathrm{Im}\, z\right|\ge R\right\} =0.
\]

\end{cor}
\noindent \textit{Proof of the lemma and its corollary.}The corollary
immediately follows from (\ref{eq:awe}), by choosing $A=X_{0}+X_{1}$
and recalling the definition of $\mathcal{F}(X_{0},X_{1})$ and the
trivial fact that $X_{j}\overset{1}{\subset}X_{0}+X_{1}$ for $j=0,1$.
Thus it remains to prove the lemma:

Let $M=\sup\left\{ \left\Vert f(z)\right\Vert :z\in\mathbb{S}\right\} $.
Given $\varepsilon>0$, we choose $R_{1}>0$ such that $\left\Vert f(j+iy)\right\Vert \le\varepsilon/e$
for $j=0,1$ and all $y\in\mathbb{R}\setminus\left(-R_{1},R_{1}\right)$.
We also choose $R_{2}$ such that $Me^{1-R_{2}^{2}}\le\varepsilon$.
Now let $z_{0}=x_{0}+iy_{0}$ be an arbitrary point of $\mathbb{S}$
for which $\mathrm{Im}\, z_{0}=y_{0}$ satisfies $\left|y_{0}\right|>R_{1}+R_{2}$
and let $\lambda\in A^{*}$ be a norm one linear functional on $A$
for which $\lambda\left(f(z_{0})\right)=\left\Vert f(z_{0})\right\Vert $.
Consider the scalar valued analytic function $g(z)=e^{\left(z-z_{0}\right)^{2}}\lambda\left(f(z)\right)$.
The point $z_{0}$ lies in the rectangle $E=\left\{ x+iy:x\in[0,1],\,\left|y-y{}_{0}\right|\le R_{2}\right\} $
and $g$ is of course continuous on $E$ and analytic in its interior.
So $\left\Vert f(z_{0})\right\Vert =\left|g(z_{0})\right|$ is dominated
by $\max\left\{ \left|g(z)\right|:z\in\partial E\right\} $, where
$\partial E$, the boundary of $E$, is the union of the ``horizonal''
set $E_{h}=\left\{ x+iy:0\le x\le1,\,\left|y-y_{0}\right|=R_{2}\right\} $
and the ``vertical'' set $E_{v}=\left\{ x+iy:x\in\left\{ 0,1\right\} ,\,\left|y-y_{0}\right|\le R_{2}\right\} $.
At every point $x+iy$ of the set $E_{h}$ we have $\left|e^{(z-z_{0})^{2}}\right|=e^{x^{2}-(y-y_{0})^{2}}\le e^{1-R_{2}^{2}}\le\varepsilon/M$,
and $\left|\lambda(f(z))\right|\le\left\Vert f(z)\right\Vert \le M$.
Therefore $\left|g(x+iy)\right|\le\varepsilon$. At every point $x+iy$
of the set $E_{v}$ we have $\left|y\right|\ge\left|y_{0}\right|-\left|y-y_{0}\right|\ge R_{1}+R_{2}-R_{2}=R_{1}$.
Therefore $\left|\lambda(f(x+iy))\right|\le\left\Vert f(x+iy)\right\Vert \le\varepsilon/e$.
Since $\left|e^{(z-z_{0})^{2}}\right|=e^{x-\left(y-y_{0}\right)^{2}}\le e$,
we obtain that $\left|g(x+iy)\right|\le\varepsilon$. We conclude
that $\left\Vert f(z_{0})\right\Vert \le\varepsilon$, which completes
the proof of the lemma. $\qed$

\smallskip{}

For any open connected set $U\subset\mathbb{C}$ we define $H^{\infty}(U)$
to be the set of all bounded analytic functions $f:U\rightarrow\mathbb{C}$
with norm $\left\Vert f\right\Vert _{H^{\infty}(U)}=\sup_{z\in U}\left|f(z)\right|$.
In fact here we are only interested in the particular spaces $H^{\infty}(\mathbb{D}^{\circ})$
and $H^{\infty}(\mathbb{S}^{\circ})$.

Of course $f\in H^{\infty}(\mathbb{S}^{\circ})$ if and only if $g=f\circ\xi$
is a function in $H^{\infty}(\mathbb{D}^{\circ})$, and the map $f\mapsto f\circ\xi$
and its inverse $g\mapsto g\circ\xi^{-1}$ define isometries between
$H^{\infty}(\mathbb{S}^{\circ})$ and $H^{\infty}(\mathbb{D}^{\circ})$.

\smallskip{}

The following theorem is a combination of well known facts about Fourier
series (see e.g., \cite{ka} Chapter 1) and the Poisson integral for
the disk (see e.g., \cite{rudin} Chapter 11). It is only a ``light''
or ``lite'' version, sufficient for our needs here%
\footnote{For our purposes here, we do not need to concern ourselves with the
deeper results which hold in this context, about the existence of
radial or even nontangential limits of the function $u$ and their
almost everywhere coincidence with the values of $h$.%
}, and we have provided an almost self contained proof. We note that
part (vi) of the theorem is closely related to Lemma \ref{3l}. 
\begin{thm}
\label{pfatou}Suppose that the measurable function $h:\mathbb{T}\rightarrow\mathbb{C}$
satisfies 
\begin{equation}
\int_{0}^{2\pi}|h(e^{it})|dt<\infty\ .\label{mybigp}
\end{equation}
 Set $\widehat{h}(n)=\frac{1}{2\pi}\int_{0}^{2\pi}e^{-int}h(e^{it})dt$
for each $n\in\mathbb{Z}$. For each $z=re^{it}$ with $r\ge0$ and
$t\in\mathbb{R}$ define 
\begin{equation}
u(z)=\sum_{n=-\infty}^{\infty}\widehat{h}(n)r^{|n|}e^{int}\label{defu}
\end{equation}
 Then:

(i) The series $\sum_{n=-\infty}^{\infty}\widehat{h}(n)r^{|n|}e^{int}$
converges absolutely for every $r\in[0,1)$ and $t\in\mathbb{R}$
and thus (\ref{defu}) defines a continuous function $u:\mathbb{D}^{\circ}\rightarrow\mathbb{C}$.

(ii) The function $u:\mathbb{D}^{\circ}\rightarrow\mathbb{C}$ is
harmonic in $\mathbb{D}^{\circ}$.

(iii) If $\widehat{h}(n)=0$ for every $n<0$ then $u$ is also analytic
in $\mathbb{D}^{\circ}$.

(iv) If $h$ is essentially bounded then 
\[
\left|u(z)\right|\le\underset{t\in[0,2\pi]}{\mathrm{ess~sup}}\,\left|h(e^{it})\right|\text{ for all }z\in\mathbb{D}^{\circ}.
\]

(v) Suppose that $f:\mathbb{D}\rightarrow\mathbb{C}$ is a continuous
function which is analytic in $\mathbb{D}^{\circ}$ and $\widehat{h}(n)=0$
for every $n<0$, and $g(e^{it})=h(e^{it})f(e^{it})$ for a.e.\ $t\in[0,2\pi]$
with Fourier coefficients $\widehat{g}(n)=\frac{1}{2\pi}\int_{0}^{2\pi}e^{-int}g(e^{it})dt$.
Then 
\[
\sum_{n=-\infty}^{\infty}\widehat{g}(n)r^{|n|}e^{int}=u(re^{it})f(re^{it})\text{ for every }r\in[0,1)\text{ and }t\in\mathbb{R},
\]
and, in particular, 
\[
\frac{1}{2\pi}\int_{0}^{2\pi}g(e^{it})dt=u(0)f(0).
\]

(vi) Let $\Gamma_{0}$ and $\Gamma_{1}$ and $\xi$ be the two arcs
and the conformal map defined in (\ref{eq:defgamj}) and (\ref{defxi})
above. Let $M_{0}$ and $M_{1}$ be positive constants. Suppose that
$\left|h(e^{it})\right|\le M_{0}$ for almost all $e^{it}\in\Gamma_{0}$
and $\left|h(e^{it})\right|\le M_{1}$ for almost all $e^{it}\in\Gamma_{1}$
and that $\widehat{h}(n)=0$ for every $n<0$. Then 
\begin{equation}
\left|u(\xi^{-1}(z))\right|\le M_{0}^{1-\mathrm{Re\,}z}M_{1}^{\mathrm{Re\,}z}\text{ for all }z\in\mathbb{S}^{\circ}.\label{flerb}
\end{equation}
 \end{thm}
\begin{rem}
In fact in our subsequent applications we will not really need the
fact (conclusion (ii) above) that $u$ is harmonic. 
\end{rem}
\textit{Proof.} Properties (i), (ii) and (iii) are rather obvious.
Property (i) is a trivial consequence of the boundedness of the sequence
$\left\{ \widehat{h}(n)\right\} _{n\in\mathbb{Z}}$. It implies in
turn that the convergence of $u_{N}(re^{it}):=\sum_{n=-N}^{N}\widehat{h}(n)r^{|n|}e^{int}$
to $u(re^{it})$ as $N$ tends to $\infty$ is uniform on every compact
subset of $\mathbb{D}^{\circ}$. Since each $u_{N}$ is harmonic,
it follows by Harnack's theorem (see e.g., \cite{rudin} p.\ 236,
Theorem 11.11) that $u$ is harmonic on $\mathbb{D}^{\circ}$. This
establishes (ii). Property (iii) follows in almost the same way since,
under the stated condition, the functions $u_{N}$ will also be analytic.

To obtain (iv) we first note that, for each fixed $r\in[0,1)$ and
$t\in\mathbb{R}$, the uniform convergence of $s\mapsto\sum_{n=-\infty}^{\infty}r^{|n|}e^{in(t-s)}$
on $\mathbb{R}$ permits us to integrate term by term in the following
calculation to obtain that

\begin{eqnarray}
\sum_{n=-\infty}^{\infty}\widehat{h}(n)r^{|n|}e^{int} & = & \sum_{n=-\infty}^{\infty}\frac{1}{2\pi}\int_{0}^{2\pi}h(e^{is})r^{|n|}e^{in(t-s)}ds\notag\\
 & = & \frac{1}{2\pi}\int_{0}^{2\pi}h(e^{is})\left[\sum_{n=-\infty}^{\infty}r^{|n|}e^{in(t-s)}\right]ds.\label{wxp}
\end{eqnarray}

\smallskip{}
 Although I have been trying obstinately to avoid dealing explicitly
with the Poisson kernel for the disk, it seems that I have no alternative
here but to carry out the (admittedly very simple) calculation (see
e.g., \cite{rudin} p.\ 111) which shows that 
\begin{equation}
\sum_{n=-\infty}^{\infty}r^{|n|}e^{in(t-s)}=\frac{1-r^{2}}{1-2r\cos(t-s)+r^{2}}.\label{eca}
\end{equation}
 The only reason that we need (\ref{eca}) is in order to be able
to affirm that 
\[
\sum_{n=-\infty}^{\infty}r^{|n|}e^{in(t-s)}\ge0\ \text{ for all \ensuremath{r\in[0,1)}and \ensuremath{s,t\in\mathbb{R}}\,.}
\]

By considering (\ref{wxp}) in the special case where $h(e^{it})=1$
for all $t$, we see that 
\[
\frac{1}{2\pi}\int_{0}^{2\pi}\left[\sum_{n=-\infty}^{\infty}r^{|n|}e^{in(t-s)}\right]ds=1.
\]
 So 
\begin{eqnarray*}
\left|u(re^{it})\right| & = & \left|\frac{1}{2\pi}\int_{0}^{2\pi}h(e^{is})\left[\sum_{n=-\infty}^{\infty}r^{|n|}e^{in(t-s)}\right]ds\right|\\
 & \le & \frac{1}{2\pi}\int_{0}^{2\pi}|h(e^{is})|\cdot\left|\sum_{n=-\infty}^{\infty}r^{|n|}e^{in(t-s)}\right|ds\\
 & \le & \underset{t\in[0,2\pi]}{\mathrm{ess~sup}}\,\left|h(e^{it})\right|\cdot\frac{1}{2\pi}\int_{0}^{2\pi}\left|\sum_{n=-\infty}^{\infty}r^{|n|}e^{in(t-s)}\right|ds=\underset{t\in[0,2\pi]}{\mathrm{ess~sup}}\,\left|h(e^{it})\right|
\end{eqnarray*}
 as required.

\smallskip{}

To prove (v) we start by rewriting the product $u(re^{it})f(re^{it})$
as a double sum 
\begin{eqnarray}
u(re^{it})f(re^{it}) & = & \sum_{n=0}^{\infty}\widehat{h}(n)r^{n}e^{int}\cdot\sum_{m=0}^{\infty}\widehat{f}(m)r^{m}e^{imt}\notag\\
 & = & \sum_{n=0}^{\infty}\left(\sum_{m=0}^{\infty}\widehat{h}(n)\widehat{f}(m)r^{m+n}e^{i(m+n)t}\right).\label{bsa}
\end{eqnarray}
 The above double sum is absolutely convergent, i.e., the following
sum 
\[
\sum_{n=0}^{\infty}\left(\sum_{m=0}^{\infty}|\widehat{h}(n)\widehat{f}(m)r^{(m+n)}e^{i(m+n)t}|\right)=\sum_{n=0}^{\infty}\left(\sum_{m=0}^{\infty}|\widehat{h}(n)||\widehat{f}(m)|r^{m+n}\right)
\]
 is finite. Therefore we can change the order of summation in (\ref{bsa})
so that in fact 
\begin{eqnarray*}
u(re^{it})f(re^{it}) & = & \sum_{k=0}^{\infty}\left(\sum_{m+n=k,m\ge0,n\ge0}\widehat{h}(n)\widehat{f}(m)r^{m+n}e^{i(m+n)t}\right)\\
 & = & \sum_{k=0}^{\infty}\left(\sum_{m+n=k,m\ge0,n\ge0}\widehat{h}(n)\widehat{f}(m)\right)r^{k}e^{ikt}.
\end{eqnarray*}

To complete the proof of part (v) we have to show that, for each $k\in\mathbb{Z}$,
\begin{equation}
\widehat{g}(k)=\left\{ \begin{array}{ccc}
{\displaystyle \sum_{m+n=k,m\ge0,n\ge0}\widehat{h}(n)\widehat{f}(m)} & , & k\ge0\\
\phantom{.}\phantom{,}\\
0 & , & k<0
\end{array}\right..\label{raqz}
\end{equation}
To do this, we first observe that $\widehat{f}(n)=\frac{1}{2\pi}\int_{0}^{2\pi}e^{-int}f(e^{it})dt=\frac{1}{2\pi i}\int_{0}^{2\pi}e^{-i(n+1)t}f(e^{it})ie^{it}dt$
is exactly the complex contour integral $\frac{1}{2\pi i}\oint_{\mathbb{T}}z^{-n-1}f(z)dz$
and so, by Cauchy's theorem, it vanishes for every negative integer
$n$. For each $k\in\mathbb{Z}$, we will calculate $\widehat{g}(k)$
with the help of the function $v:\mathbb{T}\rightarrow\mathbb{C}$
defined by $v(e^{it})=e^{ikt}\overline{f(e^{it})}$. The Fourier coefficients
of this function are $\widehat{v}(n)=\frac{1}{2\pi}\int_{0}^{2\pi}e^{-i(n-k)}\overline{f(e^{it})}dt=\overline{\widehat{f}(k-n)}$.
Both $|h|$ and $|v|$ are essentially bounded and therefore square
integrable on $\mathbb{T}$. So we can apply the generalized Parseval
identity to obtain that 
\begin{eqnarray*}
\widehat{g}(k) & = & \frac{1}{2\pi}\int_{0}^{2\pi}e^{-ikt}h(e^{it})f(e^{it})dt=\frac{1}{2\pi}\int_{0}^{2\pi}h(e^{it})\overline{v(e^{it})}dt=\sum_{n\in\mathbb{Z}}\widehat{h}(n)\overline{\widehat{v}(n)}\\
 & = & \sum_{n\in\mathbb{Z}}\widehat{h}(n)\widehat{f}(k-n).
\end{eqnarray*}
 Since the summands in this series vanish for all $n<0$ and all $n>k$,
this establishes (\ref{raqz}) for each $k\in\mathbb{Z}$ and so completes
the proof of part (v).

\smallskip{}
 For part (vi) we will once more use the particular function $\phi:\mathbb{C}\rightarrow\mathbb{C}$
which we already used in the proof of Lemma \ref{3l}. I.e., we have
$\phi(z)=\exp(\delta z^{2}-(1-z)\log M_{0}-z\log M_{1})$ where $\delta$
is a positive constant. This non-vanishing function, or rather its
restriction to $\mathbb{S}$, has all the properties required in part
(i) of Lemma \ref{lem:PhiPsi} (in the case where $A=\mathbb{C}$).
Therefore the function $f(w):=\phi(\xi(w))$ is analytic and non zero
at every point $w\in\mathbb{D}\backslash\{1,e^{i2\pi\theta}\}$ and
extends to a continuous function on $\mathbb{D}$ with $f(1)=f(e^{i2\pi\theta})=0$.
I.e., it satisfies the hypotheses of (v). 

Furthermore (recalling (\ref{idrvba})) we see that, for $j=0,1$,
\begin{equation}
\sup_{e^{it}\in\Gamma_{j}}|f(e^{it})|=\sup_{t\in\mathbb{R}}|\phi(j+it)|=e^{\delta j}/M_{j}\text{ }.\label{bdietma}
\end{equation}

\smallskip{}
 It follows from (\ref{bdietma}) and our hypotheses on $h$, that
the function $g(e^{it}):=h(e^{it})f(e^{it})$ satisfies $\mathrm{ess~sup}_{t\in[0,2\pi]}|g(e^{it})|\le e^{\delta}$.
Thus, applying part (iv), but to the function $g$ instead of $h$,
and then applying part (v), we see that 
\[
e^{\delta}\ge\left|\sum_{n=-\infty}^{\infty}\widehat{g}(n)r^{|n|}e^{int}\right|=|u(re^{it})f(re^{it})|\text{ for every }r\in[0,1)\text{ and }t\in\mathbb{R}.
\]

Given any $z\in\mathbb{S}^{\circ}$, choose $r\in[0,1)$ and $t\in\mathbb{R}$
such that $z=\xi(re^{it})$. Then 
\[
\left|u(\xi^{-1}(z))\right|=\left|u(re^{it})\right|\le\frac{e^{\delta}}{\left|f(re^{it})\right|}=\frac{e^{\delta}}{\left|\phi(z)\right|}=\frac{e^{\delta}}{|e^{\delta z^{2}}|M_{0}^{-1+\mathrm{Re\,}z}M_{1}^{-\mathrm{Re\,}z}}.
\]
 Since we can choose the positive number $\delta$ to be as small
as we please, these estimates immediately imply (\ref{flerb}) and
complete the proof of (vi) and thus of Theorem \ref{pfatou}. \qed

\smallskip{}

We will also need another result about functions in $H^{\infty}(\mathbb{S}^{\circ})$:

\smallskip{}

\begin{lem}
\label{bapf}Let $h$ be a function in $H^{\infty}(\mathbb{S}^{\circ})$.
Then there exists a function $\phi:\mathbb{S}\rightarrow\mathbb{C}$
with the following properties:

(i) $\left|\phi(z_{1})-\phi(z_{2})\right|\le\left|z_{1}-z_{2}\right|\cdot\sup_{\zeta\in\mathbb{S}^{\circ}}|h(\zeta)|$
for all $z_{1},z_{2}\in\mathbb{S}.$

(ii) $\phi$ is analytic in $\mathbb{S}^{\circ}$ and satisfies $\phi^{\prime}(z)=h(z)$
for each $z\in\mathbb{S}^{\circ}$,

(iii) The estimate 
\[
\left|\phi(j+it_{1})-\phi(j+it_{2})\right|\le\left|t_{1}-t_{2}\right|\cdot\lim_{r\searrow0}\left(\sup\left\{ \left|h(z)\right|:z\in\mathbb{S}\,,\text{ }0<\left|\mathrm{Re\,}z-j\right|<r\right\} \right)
\]
 holds for $j=0,1$ and all $t_{1},t_{2}\in\mathbb{R}$.

(iv) $\left|\phi(z)\right|\le\left|z-\frac{1}{2}\right|\cdot\sup_{\zeta\in\mathbb{S}^{\circ}}|h(\zeta)|$
for all $z\in\mathbb{S}$. 
\end{lem}
\textit{Proof.} We define $\phi(z)$ for all $z\in\mathbb{S}$ by
the two formul\ae 
\begin{equation}
\phi(z)=\int_{1/2}^{z}h(\zeta)d\zeta\text{ for each }z\in\mathbb{S}^{\circ}\label{rty}
\end{equation}
 and 
\begin{equation}
\phi(j+it)=\lim_{n\rightarrow\infty}\int_{1/2}^{j+it+(1/2-j)/n}h(\zeta)d\zeta\text{ for each }j=0,1\text{ and }t\in\mathbb{R}.\label{zrty}
\end{equation}
 In (\ref{rty}) the contour integration is performed along a contour,
any contour, in $\mathbb{S}^{\circ}$ from $1/2$ to $z$. This ensures
that property (ii) holds. We obtain (i), at least in the case where
$z_{1}$ and $z_{2}$ are both in $\mathbb{S}^{\circ}$, by choosing
a contour from $1/2$ to $z_{1}$ which includes a straight line segment
from $z_{2}$ to $z_{1}$. This shows that $\mathbb{\phi}$ is uniformly
continuous on $\mathbb{S}^{\circ}$ and therefore we can extend our
definition of $\phi$ to all of $\mathbb{S}$ in a suitable way. This
will of course be the unique continuous extension of $\phi$ to $\mathbb{S}$
whose values on $\partial\mathbb{S}$ must necessarily equal those
given by the formula (\ref{zrty}). It is now clear that $\phi$ satisfies
(i) on all of $\mathbb{S}$, and property (iv) follows obviously from
property (i) together with the fact that $\phi(1/2)=0$.

It remains to check that property (iii) also holds: For $j=0,1$ and
for each real $t_{1}$ and $t_{2}$ we have 
\begin{eqnarray}
\phi(j+it_{2})-\phi(j+it_{1}) & = & \lim_{n\rightarrow\infty}\phi\left(j+it_{2}+\frac{1/2-j}{n}\right)-\phi\left(j+it_{1}+\frac{1/2-j}{n}\right)\notag\\
 & = & i\lim_{n\rightarrow\infty}\int_{t_{1}}^{t_{2}}h\left(j+\frac{1/2-j}{n}+it\right)dt\text{ }.\label{asl}
\end{eqnarray}

\smallskip{}
 Given $\epsilon>0$, we can choose $r_{\epsilon}>0$ such that 
\begin{eqnarray*}
 &  & \sup\left\{ \left|h(z)\right|:z\in\mathbb{S},\text{ }0<\left|\mathrm{Re\,}z-j\right|<r_{\epsilon}\right\} \\
 & \le & \epsilon+\lim_{r\searrow0}\left(\sup\left\{ \left|h(z)\right|:z\in\mathbb{S},\text{ }0<\left|\mathrm{Re\,}z-j\right|<r\right\} \right).
\end{eqnarray*}
 So, for all $n>1/r_{\epsilon}$, we have 
\[
\left|h(j+(1/2-j)/n+it)\right|\le\epsilon+\lim_{r\searrow0}\left(\sup\left\{ \left|h(z)\right|:z\in\mathbb{S},\text{ }0<\left|\mathrm{Re\,}z-j\right|<r\right\} \right).
\]
 This, combined with (\ref{asl}), shows that 
\[
\left|\phi(j+it_{2})-\phi(j+it_{1})\right|\le|t_{1}-t_{2}|\left(\epsilon+\lim_{r\searrow0}\left(\sup\left\{ \left|h(z)\right|:z\in\mathbb{S},\text{ }0<\left|\mathrm{Re\,}z-j\right|<r\right\} \right)\right).
\]
 Since $\epsilon$ can be chosen arbitrarily small, we obtain property
(iii). \qed

\smallskip{}

\subsection{\protect \label{altdeflower}An alternative definition of the space
$[X_{0},X_{1}]_{\theta}$, and another preliminary result.}
\begin{defn}
\label{deffd} We define $\mathcal{F}_{\mathbb{D}}(X_{0},X_{1})$
to be the space of all continuous functions $f:\mathbb{D}\rightarrow X_{0}+X_{1}$
which are analytic in $\mathbb{D}^{\circ}$ and satisfy $f(1)=f(e^{i2\pi\theta})=0$
and such that, for $j=0,1$, the restriction of $f$ to the closed
arc $\overline{\Gamma_{j}}$ is a continuous $X_{j}$ valued function.
We norm this space by 
\begin{equation}
\left\Vert f\right\Vert _{\mathcal{F}_{\mathbb{D}}(X_{0},X_{1})}=\max\left\{ \max_{z\in\Gamma_{0}}\left\Vert f(z)\right\Vert _{X_{0}},\max_{z\in\Gamma_{1}}\left\Vert f(z)\right\Vert _{X_{1}}\right\} .\label{tulfamw}
\end{equation}

\end{defn}
We claim that any given function $g:\mathbb{S}\rightarrow X_{0}+X_{1}$
is an element of $\mathcal{F}(X_{0},X_{1})$ if and only if the function
$g\circ\xi:\mathbb{D}\backslash\{1,e^{i2\pi\theta}\}\rightarrow X_{0}+X_{1}$,
can be extended continuously to all of $\mathbb{D}$ and this extension
is an element of $\mathcal{F}_{\mathbb{D}}(X_{0},X_{1})$. (Here $\xi$
is of course once more the function defined by (\ref{defxi}).) In
that case we also have 
\begin{equation}
\left\Vert g\right\Vert _{\mathcal{F}(X_{0},X_{1})}=\left\Vert g\circ\xi\right\Vert _{\mathcal{F}_{\mathbb{D}}(X_{0},X_{1})}.\label{eq:IsomSD}
\end{equation}
One ingredient of the straightforward proof of this claim is Corollary
\ref{cor:VanishAtInfty}, and others are (\ref{idrvba}), (\ref{eq:upstrip})
and (\ref{eq:DownStrip}). The remaining ingredients are obvious arguments
very similar to those needed to prove Lemma \ref{lem:PhiPsi}. %
\begin{comment}
\textit{\textcolor{red}{Maybe be more explicit about this in future
versions.}}
\end{comment}

Since $\xi(0)=\theta$, we deduce immediately that 
\begin{equation}
[X_{0},X_{1}]_{\theta}=\left\{ f(0):f\in\mathcal{F}_{\mathbb{D}}(X_{0},X_{1})\right\} \label{rtza}
\end{equation}
and 
\begin{equation}
\left\Vert x\right\Vert _{[X_{0},X_{1}]_{\theta}}=\inf\left\{ \left\Vert f\right\Vert _{\mathcal{F}_{\mathbb{D}}(X_{0},X_{1})}:f\in\mathcal{F}_{\mathbb{D}}(X_{0},X_{1}),f(0)=x\right\} .\label{rtzaa}
\end{equation}

Calderón's space $\mathcal{G}(X_{0},X_{1})$, defined on p.\ 116
of \cite{ca}, which (cf.~Definition \ref{def:3DifferentGs}) we
denote here by $\mathcal{G}_{\mathbf{Ca}}\left(X_{0},X_{1}\right)$,
is a very useful dense subspace of $\mathcal{F}(X_{0},X_{1})$. We
will need a more or less analogous dense subspace of $\mathcal{F}_{\mathbb{D}}(X_{0},X_{1})$
which in fact will correspond exactly to Stafney's space $\mathcal{G}(X_{0},X_{1})$,
defined on p.\ 334 of \cite{stafney} and which we denote here by
$\mathcal{G}_{\mathbf{St}}\left(X_{0},X_{1}\right)$. 
\begin{defn}
\label{badefg}Let $\mathcal{G}_{\mathbb{D}}(X_{0},X_{1})$ be the
subspace of $\mathcal{F}_{\mathbb{D}}(X_{0},X_{1})$ consisting of
all elements $f:\mathbb{D}\rightarrow X_{0}\cap X_{1}$ which are
finite sums of the form 
\begin{equation}
f(w)=\sum_{k=1}^{N}\psi_{k}(w)a_{k}\label{defg}
\end{equation}
where each $a_{k}\in X_{0}\cap X_{1}$ and each $\psi_{k}:\mathbb{D}\rightarrow\mathbb{C}$
is a continuous function on $\mathbb{D}$ which is analytic in $\mathbb{D}^{\circ}$
and satisfies $\psi_{k}(1)=\psi_{k}(e^{i2\pi\theta})=0$. 
\end{defn}
For each $\psi_{k}$ in (\ref{defg}), let us define the function
$\phi_{k}:\mathbb{S}\to\mathbb{C}$ by $\phi_{k}(z)=\psi_{k}(\xi^{-1}(z))$
where $\xi^{-1}$ is the map of $\mathbb{S}$ onto $\mathbb{D}\setminus\left\{ 1,e^{i\pi\theta}\right\} $
defined above by (\ref{eq:defxiinv}). By part (ii) of Lemma \ref{lem:PhiPsi},
in the case $A=\mathbb{C}$, $\psi_{k}$ has the same properties that
are demanded of the functions $c_{n}$ in the definition near the
bottom of p.~334 of \cite{stafney}. It follows that, for each $f\in\mathcal{G}_{\mathbb{D}}(X_{0},X_{1})$,
the function $g=f\circ\xi^{-1}$ is an element of Stafney's space
$\mathcal{G}_{\mathbf{St}}\left(X_{0},X_{1}\right)$. Conversely,
for every function $g\in\mathcal{G}_{\mathbf{St}}\left(X_{0},X_{1}\right)$,
we can apply part (i) of Lemma \ref{lem:PhiPsi} to the scalar functions
in the formula for $g$ to deduce that the function $f$, defined
by 
\begin{equation}
f(w)=\left\{ \begin{array}{lll}
g(\xi(w)) & , & w\in\mathbb{D}\backslash\{1,e^{i2\pi\theta}\}\\
0 & , & w=1,e^{i2\pi\theta}
\end{array}\right.,\label{eq:aqwq}
\end{equation}
is an element of $\mathcal{G}_{\mathbb{D}}(X_{0},X_{1})$.

Definition \ref{badefg} gives a somewhat larger class than would
be obtained by taking the image of the space $\mathcal{G}_{\mathbf{Ca}}(X_{0},X_{1})$
(as defined via Definition \ref{def:3DifferentGs} and \cite[p. 116]{ca})
under composition with the function $\xi$. In fact the class of functions
$f:\mathbb{D}\to X_{0}\cap X_{1}$ obtained via (\ref{eq:aqwq}) for
all functions $g:\mathbb{S}\rightarrow X_{0}\cap X_{1}$ which are
elements of $\mathcal{G}_{\mathbf{Ca}}(X_{0},X_{1})$, is precisely
the subspace of $\mathcal{G}_{\mathbb{D}}(X_{0},X_{1})$ where the
functions $\psi_{k}$ in (\ref{defg}) are each required to be of
the special form $\psi_{k}(w)=e^{\delta_{k}\xi(w)^{2}+i\lambda_{k}\xi(w)}$
for some constants $\delta_{k}>0$ and $\lambda_{k}\in\mathbb{R}$
for all $w\in\mathbb{D}\backslash\{1,e^{i2\pi\theta}\}$. This class
is dense in $\mathcal{F}_{\mathbb{D}}(X_{0},X_{1})$ because of (\ref{eq:IsomSD})
and the density (see \cite{ca} Section 9.2 p.\ 116 and Section 29.2
pp.~132--133) of $\mathcal{G}_{\mathbf{Ca}}(X_{0},X_{1})$ in $\mathcal{F}(X_{0},X_{1})$.
So of course $\mathcal{G}_{\mathbb{D}}(X_{0},X_{1})$ must also be
dense in $\mathcal{F}_{\mathbb{D}}(X_{0},X_{1})$.

\begin{rem}
Our definitions here of $\mathcal{F}_{\mathbb{D}}(X_{0},X_{1})$ and
$\mathcal{G}_{\mathbb{D}}(X_{0},X_{1})$ depend of course on our choice
of the parameter $\theta\in(0,1)$. In some contexts, for example
where, unlike here, it is necessary to simultaneously consider different
values of $\theta$, it might be more natural to first define $\xi$
and $\mathcal{F}_{\mathbb{D}}(X_{0},X_{1})$ for the particular choice
$\theta=1/2$. In that case, instead of (\ref{rtza}) and (\ref{rtzaa})
we would have $[X_{0},X_{1}]_{\alpha}=\left\{ f(\xi^{-1}(\alpha)):f\in\mathcal{F}_{\mathbb{D}}(X_{0},X_{1})\right\} $
and 
\[
\left\Vert x\right\Vert _{[X_{0},X_{1}]_{\alpha}}=\inf\left\{ \left\Vert f\right\Vert _{\mathcal{F}_{\mathbb{D}}(X_{0},X_{1})}:f\in\mathcal{F}_{\mathbb{D}}(X_{0},X_{1}),f(\xi^{-1}(\alpha))=x\right\} 
\]
 for each $\alpha\in(0,1)$. 
\end{rem}
The following result is essentially equivalent to the third inequality
stated in Section 9.4 on p.\ 117 of \cite{ca}, and proved in Section
29.4 on pp.\ 134--135 of \cite{ca}. But we offer the reader a self
contained and perhaps in some ways simpler%
\footnote{One thing which might be considered a simplification here is that,
in contrast to \cite{ca} and \cite{ccrsw}, we only have to consider
the obvious and explicit harmonic conjugates of harmonic polynomials
instead of harmonic conjugates of more general smooth functions.%
} proof using almost the same approach as in \cite{ca}. (Cf.\ also
\cite{ccrsw} Proposition 2.4 pp.~209--210).

\smallskip{}

\begin{lem}
\label{gie}The inequality 
\begin{equation}
\left\Vert f(0)\right\Vert _{[X_{0},X_{1}]_{\theta}}\le\frac{1}{2\pi}\left(\int_{2\pi\theta}^{2\pi}\left\Vert f(e^{it})\right\Vert _{X_{0}}dt+\int_{0}^{2\pi\theta}\left\Vert f(e^{it})\right\Vert _{X_{1}}dt\right)\label{aqz}
\end{equation}

holds for every $f\in\mathcal{F}_{\mathbb{D}}(X_{0},X_{1})$. 
\end{lem}
\textit{Proof.} It will be notationally convenient to introduce the
continuous function $h:\mathbb{T}\rightarrow[0,\infty)$ defined by
\begin{equation}
h(e^{it})=\left\{ \begin{array}{lll}
\left\Vert f(e^{it})\right\Vert _{X_{0}} & , & 2\pi\theta<t<2\pi\\
\left\Vert f(e^{it})\right\Vert _{X_{1}} & , & 0<t<2\pi\theta\\
0 & , & t=0,\text{ }2\pi\theta.
\end{array}\right.\label{werp}
\end{equation}
 and so to rewrite (\ref{aqz}) as 
\begin{equation}
\left\Vert f(0)\right\Vert _{[X_{0},X_{1}]_{\theta}}\le\frac{1}{2\pi}\int_{0}^{2\pi}h(e^{it})dt.\label{dyam}
\end{equation}

By Jensen's inequality (\cite{rudin} Theorem 3.3 p.\ 62) and the
convexity of the exponential function, we have 
\begin{equation}
\exp\left(\frac{1}{2\pi}\int_{0}^{2\pi}\log\left(h(e^{it})+\rho\right)dt\right)\le\frac{1}{2\pi}\left(\int_{0}^{2\pi}\left(h(e^{it})+\rho\right)dt\right)\label{npfe}
\end{equation}
 for each positive number $\rho$. So our main step will be to prove
that 
\begin{equation}
\left\Vert f(0)\right\Vert _{[X_{0},X_{1}]_{\theta}}\le\exp\left(\frac{1}{2\pi}\int_{0}^{2\pi}\log\left(h(e^{it})+\rho\right)dt\right)\label{raq}
\end{equation}
 for each $\rho>0$. Then we can obtain (\ref{dyam}) by simply substituting
in (\ref{npfe}), and then taking the limit as $\rho$ tends to $0$.

Choose an arbitrary positive number $\epsilon$. Since $e^{it}\mapsto\log\left(h(e^{it})+\rho\right)$
is a continuous real valued function on $\mathbb{T}$, there exists
(cf.\ e.g., \cite{ka} p.\ 15) a real valued trigonometric polynomial
$p(e^{it})=\sum_{n=-N}^{N}c_{n}e^{int}$ such that 
\begin{equation}
\log\left(h(e^{it})+\rho\right)\le p(e^{it})+\epsilon\le\log\left(h(e^{it})+\rho\right)+2\epsilon\label{rimlo}
\end{equation}
for all $e^{it}\in\mathbb{T}$. The fact that our polynomial is real
valued implies that $\overline{c_{n}}=c_{-n}$ for each $n$. A trivial
calculation using this last formula shows that the functions $u:\mathbb{C}\rightarrow\mathbb{C}$
and $v:\mathbb{C}\rightarrow\mathbb{C}$ defined by 
\[
u(re^{it})=\sum_{n=-N}^{N}c_{n}r^{|n|}e^{int}\text{ and }v(re^{it})=\frac{1}{i}\sum_{n=1}^{N}c_{n}r^{|n|}e^{int}-\frac{1}{i}\sum_{n=-N}^{-1}c_{n}r^{|n|}e^{int}\text{ }
\]
for all $r\ge0$ and $t\in\mathbb{R}$ are both real valued. It is
also clear that the function $w:\mathbb{C}\rightarrow\mathbb{C}$
defined by $w(z)=u(z)+iv(z)$ is a polynomial in $z$ and thus analytic
for all $z$.

This last fact ensures that the function $g(z):=$ $e^{-w(z)-\epsilon}f(z)$
is an element of $\mathcal{F}_{\mathbb{D}}(X_{0},X_{1})$. Furthermore,
for each $e^{it}\in\mathbb{T}$ we have 
\[
\left|e^{-w(z)-\epsilon}\right|=e^{-u(e^{it})-\epsilon}=e^{-p(e^{it})-\epsilon}\le e^{-\log(h(e^{it})+\rho)}=\frac{1}{h(e^{it})+\rho}.
\]
Thus, for $t\in(2\pi\theta,2\pi)$, i.e., for $z=e^{it}\in\Gamma_{0}$,
we have, recalling (\ref{werp}), that 
\[
\left\Vert g(z)\right\Vert _{X_{0}}=\left|e^{-w(z)-\epsilon}\right|\left\Vert f(e^{it})\right\Vert _{X_{0}}\le\frac{\left\Vert f(e^{it})\right\Vert _{X_{0}}}{\left\Vert f(e^{it})\right\Vert _{X_{0}}+\rho}\le1.
\]
 Exactly analogous considerations show that $\left\Vert g(e^{it})\right\Vert _{X_{1}}\le1$
for each $e^{it}\in\Gamma_{1}$ and so we have $\left\Vert g(0)\right\Vert _{[X_{0},X_{1}]_{\theta}}\le\left\Vert g\right\Vert _{\mathcal{F}_{\mathbb{D}}(X_{0},X_{1})}\le1$.
This means that 
\begin{equation}
\left\Vert f(0)\right\Vert _{[X_{0},X_{1}]_{\theta}}=\left|e^{w(0)+\epsilon}\right|\left\Vert g(0)\right\Vert _{[X_{0},X_{1}]_{\theta}}\le\left|e^{w(0)+\epsilon}\right|=e^{u(0)+\epsilon}.\label{ubfd}
\end{equation}
Since $u$ is a harmonic function, (or, even more simply, since $u(0)=c_{0}$
and $\int_{0}^{2\pi}e^{int}dt=0$ for every non zero $n\in\mathbb{Z}$)
we have $u(0)=\frac{1}{2\pi}\int_{0}^{2\pi}u(e^{it})dt=\frac{1}{2\pi}\int_{0}^{2\pi}p(e^{it})dt$.
This, combined with (\ref{ubfd}) and (\ref{rimlo}) shows that 
\[
\left\Vert f(0)\right\Vert _{[X_{0},X_{1}]_{\theta}}\le\exp\left(\frac{1}{2\pi}\int_{0}^{2\pi}(\log(h(e^{it})+\rho)+\epsilon)dt+\epsilon\right).
\]
 Since $\epsilon$ can be chosen arbitrarily small, this gives us
(\ref{raq}) and so completes the proof of the lemma. \qed

\smallskip{}

We conclude this subsection with a lemma which in fact is not required
for the proof of Theorem \ref{crt}. We will probably need it, or
rather its almost immediate corollary, in a paper which is currently
in preparation. The lemma is a very special case, but one sufficient
for our purposes in that paper, of a quite standard result, namely,
Theorem 15.19 of \cite{rudin} p.\ 308. We waited till now, rather
than including it earlier in Subsection \ref{svaf}, since its proof
can be almost immediately deduced from part of the proof of Lemma
\ref{gie}.
\begin{lem}
\label{postgie}Suppose that the function $\psi:\mathbb{D}\rightarrow\mathbb{C}$
is continuous in $\mathbb{D}$ and is analytic in $\mathbb{D}^{\circ}$.
Suppose further that $\psi(1)=\psi(e^{i2\pi\theta})=0$. Then 
\begin{equation}
\left|\psi(0)\right|\le\exp\left(\frac{1}{2\pi}\int_{0}^{2\pi}\log\left|\psi(e^{it})\right|dt\right).\label{rffr}
\end{equation}
 \end{lem}
\begin{cor}
\label{cgie}Every function $\psi$ satisfying the conditions of Lemma
\ref{postgie} also satisfies 
\begin{equation}
\left|\psi(0)\right|\le\left(\frac{1}{2\pi(1-\theta)}\int_{2\pi\theta}^{2\pi}\left|\psi(e^{it})\right|dt\right)^{1-\theta}\left(\frac{1}{2\pi\theta}\int_{0}^{2\pi\theta}\left|\psi(e^{it})\right|dt\right)^{\theta}.\label{udam}
\end{equation}
 
\end{cor}
\textit{Proof of Lemma \ref{postgie}.} Let $X$ be a non trivial
Banach space. (We can for example suppose that $X=\mathbb{C}$.) Suppose
that $X_{0}=X_{1}=X$ with equality of norms. Then it is a very simple
exercise to show that $[X_{0},X_{1}]_{\theta}=X$ with equality of
norms for each $\theta\in(0,1)$. Let $a$ be an element of $X$ with
$\left\Vert a\right\Vert _{X}=1$. Then the function $f(z)=\psi(z)a$
is in $\mathcal{F}_{\mathbb{D}}(X_{0},X_{1})=\mathcal{F}_{\mathbb{D}}(X,X)$.
In this case the function $h$ defined by (\ref{werp}) satisfies
\begin{equation}
h(e^{it})=\left|\psi(e^{it})\right|\text{ for all }t\in[0,2\pi]\text{ .}\label{cmejct}
\end{equation}
 for all $t\in[0,2\pi]$. Furthermore 
\begin{equation}
\left\Vert f(0)\right\Vert _{[X_{0},X_{1}]_{\theta}}=\left\Vert f(0)\right\Vert _{X}=\left|\psi(0)\right|.\label{lywc}
\end{equation}
To obtain (\ref{rffr}), all we have to do now is to substitute (\ref{cmejct})
and (\ref{lywc}) in (\ref{raq}) and let $\rho$ tend to $0$. \qed

\textit{Proof of Corollary \ref{cgie}.} We write 
\begin{eqnarray*}
 &  & \exp\left(\frac{1}{2\pi}\int_{0}^{2\pi}\log\left|\psi(e^{it})\right|dt\right)\\
 & = & \left[\exp\left(\frac{1}{2\pi}\int_{2\pi\theta}^{2\pi}\log\left|\psi(e^{it})\right|dt\right)\right]\cdot\left[\exp\left(\frac{1}{2\pi}\int_{0}^{2\pi\theta}\log\left|\psi(e^{it})\right|dt\right)\right]\\
 & = & \left[\exp\left(\frac{1}{2\pi(1-\theta)}\int_{2\pi\theta}^{2\pi}\log\left|\psi(e^{it})\right|dt\right)\right]^{1-\theta}\cdot\left[\exp\left(\frac{1}{2\pi\theta}\int_{0}^{2\pi\theta}\log\left|\psi(e^{it})\right|dt\right)\right]^{\theta}
\end{eqnarray*}
 and then simply apply Jensen's inequality to each of the two factors
in the last line. \qed

\smallskip{}

\begin{rem}
We have already indicated that (\ref{rffr}) (and therefore also (\ref{udam}))
also holds for much more general functions $\psi$. In particular
the condition $\psi(1)=\psi(e^{i2\pi\theta})=0$ is completely unnecessary
and it is almost embarrassing to impose it. However the proof of (\ref{rffr})
without imposing this condition would not have the great convenience
of using the proof of Lemma \ref{gie}. 
\end{rem}

\subsection{\protect An alternative definition of the Calderón ``upper'' method.}
\begin{defn}
\label{upper}Let $(X_{0},X_{1})$ be a regular Banach couple. Let
$\mathcal{H}(X_{0}^{\#},X_{1}^{\#})$ be the space of functions $h:\mathbb{S}^{\circ}\rightarrow X_{0}^{\#}+X_{1}^{\#}$
with the following properties:

(i) For each fixed $b\in X_{0}\cap X_{1}$, the function $z\mapsto\left\langle b,h(z)\right\rangle $
is an element of $H^{\infty}(\mathbb{S}^{\circ})$.

(ii) There exists a constant $C>0$ such that 
\begin{equation}
\left|\left\langle b,h(z)\right\rangle \right|\le C\left\Vert b\right\Vert _{X_{0}}^{1-\mathrm{Re\,}z}\left\Vert b\right\Vert _{X_{1}}^{\mathrm{Re\,}z}\text{ for all }z\in\mathbb{S}^{\circ}\text{ and for each }b\in X_{0}\cap X_{1}.\label{baconv}
\end{equation}

For each $h\in\mathcal{H}(X_{0}^{\#},X_{1}^{\#})$, the quantity $\left\Vert h\right\Vert _{\mathcal{H}(X_{0}^{\#},X_{1}^{\#})}$
is defined by 
\begin{equation}
\left\Vert h\right\Vert _{\mathcal{H}(X_{0}^{\#},X_{1}^{\#})}=\inf C\label{Aq}
\end{equation}
 where the infimum is taken over all constants $C$ which satisfy
condition (ii). \end{defn}
\begin{rem}
\label{newntl} Note that (\ref{baconv}) implies that, for $j=0,1$,
\[
\lim_{r\searrow0}\sup\left\{ \left|\left\langle b,h(z)\right\rangle \right|:z\in\mathbb{S}^{\circ}\,,\left|\mathrm{Re\,}z-j\right|<r\right\} \le C\left\Vert b\right\Vert _{X_{j}}.
\]

\end{rem}

\begin{rem}
\label{rem:UseHilPhil}In view of Fact \ref{sumq}, $X_{0}\cap X_{1}$
is a determining manifold in $\left(X_{0}^{\#}+X_{1}^{\#}\right)^{*}$
(cf.~Definition \ref{def:DeterminingManifold}) and so condition
(i) of Definition \ref{upper} guarantees that the function $h$ is
an analytic $X_{0}^{\#}+X_{1}^{\#}$ valued function on $\mathbb{S}^{\circ}$
(cf.~Definition \ref{def:HolomVecValued}). Therefore, by ``classical''
results which are recalled in Subsection \ref{hilphi}, $h$ is also
analytic in various ways which are ``stronger'' than the way specified
in condition (i). (Cf.~also Remark \ref{rem:WeakerAnalyticity}.)
\end{rem}

\begin{rem}
\label{2newntl} It might in some sense have seemed more natural in
Definition \ref{upper} to impose, instead of (ii), a condition on
the boundary values of $h$, i.e., to require that $h(j+it)$ is defined
in some perhaps ``weak'' sense for a.e.\ $t$ and is an essentially
bounded $X_{j}^{\#}$ valued function. Indeed, in the more general
context of \cite{ccrsw} (see \cite{ccrsw} Definition 3.1 pp.\ 211--212)
this is apparently the most natural way to proceed. But for our purposes
here it turns out to be technically simpler to impose the (essentially
equivalent) condition (ii) and thus to simply bypass the issues of
the existence of weak star nontangential limits etc. 
\end{rem}
\smallskip{}

It can be verified that the functional defined by (\ref{Aq}) is a
norm on the space $\mathcal{H}(X_{0}^{\#},X_{1}^{\#})$ and that $\mathcal{H}(X_{0}^{\#},X_{1}^{\#})$
is a Banach space with respect to this norm. But these facts will
also follow in a moment from Theorem \ref{babaeq}. Let us first note
a connection with the space $\mathcal{F}\left(X_{0}^{\#},X_{1}^{\#}\right)$. 
\begin{fact}
\label{fac:cbfh} $\mathcal{F}\left(X_{0}^{\#},X_{1}^{\#}\right)\subset\mathcal{H}\left(X_{0}^{\#},X_{1}^{\#}\right)$
and 
\begin{equation}
\left\Vert f\right\Vert _{\mathcal{H}\left(X_{0}^{\#},X_{1}^{\#}\right)}\le\left\Vert f\right\Vert _{\mathcal{F}\left(X_{0}^{\#},X_{1}^{\#}\right)}\,\mbox{for each }f\in\mathcal{F}\left(X_{0}^{\#},X_{1}^{\#}\right).\label{eq:hlf}
\end{equation}

\end{fact}
\textit{Proof.} It follows immediately from the definition of $\mathcal{F}\left(X_{0}^{\#},X_{1}^{\#}\right)$
that, for each $b\in X_{0}\cap X_{1}$ and each $f\in\mathcal{F}\left(X_{0}^{\#},X_{1}^{\#}\right)$,
the function $z\mapsto\left\langle b,f(z)\right\rangle $ is continuous
and bounded on $\mathbb{S}$ and analytic in $\mathbb{S}^{\circ}$
and satisfies 
\[
\left|\left\langle b,f(j+it)\right\rangle \right|\le\left\Vert b\right\Vert _{X_{j}}\left\Vert f(j+it)\right\Vert _{X_{j}^{\#}}\le\left\Vert b\right\Vert _{X_{j}}\left\Vert f\right\Vert _{\mathcal{F}\left(X_{0}^{\#},X_{1}^{\#}\right)}\,\mbox{for }j=0,1\mbox{ \mbox{and all}}t\in\mathbb{R}.
\]
So part (i) of Lemma \ref{3l} gives us that 
\[
\left|\left\langle b,f(z)\right\rangle \right|\le\left\Vert f\right\Vert _{\mathcal{F}\left(X_{0}^{\#},X_{1}^{\#}\right)}\left\Vert b\right\Vert _{B_{0}}^{1-\mathrm{Re\,}z}\left\Vert b\right\Vert _{B_{1}}^{\mathrm{Re\,}z}\text{ for all }z\in\mathbb{S}^{\circ},
\]
which completes the proof. $\qed$

\smallskip{}

\begin{thm}
\label{babaeq}Let $(X_{0},X_{1})$ be a regular Banach couple.

(i) Suppose that the function $f:\mathbb{S}\rightarrow X_{0}^{\#}+X_{1}^{\#}$
is an element of $\overline{\mathcal{F}}(X_{0}^{\#},X_{1}^{\#})$.
Then the derivative $\left(\left.f\right|_{\mathbb{S}^{\circ}}\right)^{\prime}$
of its restriction $\left.f\right|_{\mathbb{S}^{\circ}}$ to $\mathbb{S}^{\circ}$
is an element of $\mathcal{H}(X_{0}^{\#},X_{1}^{\#})$ and satisfies
\begin{equation}
\left\Vert f\right\Vert _{\overline{\mathcal{F}}(X_{0}^{\#},X_{1}^{\#})}=\left\Vert \left(\left.f\right|_{\mathbb{S}^{\circ}}\right)^{\prime}\right\Vert _{\mathcal{H}(X_{0}^{\#},X_{1}^{\#})}.\label{wwx}
\end{equation}

(ii) Conversely, suppose that the function $h:\mathbb{S}^{\circ}\rightarrow X_{0}^{\#}+X_{1}^{\#}$
is an element of $\mathcal{H}(X_{0}^{\#},X_{1}^{\#})$. Then there
exists a function $f:\mathbb{S}\rightarrow X_{0}^{\#}+X_{1}^{\#}$
in $\overline{\mathcal{F}}(X_{0}^{\#},X_{1}^{\#})$ whose restriction
$\left.f\right|_{\mathbb{S}^{\circ}}$ to $\mathbb{S}^{\circ}$ has
its derivative equal to $h$ and 
\begin{equation}
\left\Vert f\right\Vert _{\overline{\mathcal{F}}(X_{0}^{\#},X_{1}^{\#})}=\left\Vert h\right\Vert _{\mathcal{H}(X_{0}^{\#},X_{1}^{\#})}.\label{dlo}
\end{equation}
 
\end{thm}
\smallskip{}

\smallskip{}
 An immediate and obvious consequence of this theorem is
\begin{cor}
\label{corbab}The space $[X_{0}^{\#},X_{1}^{\#}]^{\theta}$ coincides
with the set of all elements $y\in X_{0}^{\#}+X_{1}^{\#}$ which arise
as the values $y=h(\theta)$ at $\theta$ of some element $h\in\mathcal{H}(X_{0}^{\#},X_{1}^{\#})$.
Furthermore 
\begin{equation}
\left\Vert y\right\Vert _{[X_{0}^{\#},X_{1}^{\#}]^{\theta}}=\inf\left\{ \left\Vert h\right\Vert _{\mathcal{H}(X_{0}^{\#},X_{1}^{\#})}:h\in\mathcal{H}(X_{0}^{\#},X_{1}^{\#}),h(\theta)=y\right\} .\label{eq:wfnt}
\end{equation}

\end{cor}

This in turn, in combination with Fact \ref{fac:cbfh}, immediately
implies the norm one inclusion 
\begin{equation}
[X_{0}^{\#},X_{1}^{\#}]_{\theta}\overset{1}{\subset}[X_{0}^{\#},X_{1}^{\#}]^{\theta}.\label{eq:Norm1Inclusion}
\end{equation}
In fact this inclusion also holds for arbitrary couples, i.e, those
which do not necessarily consist of dual spaces. This is stated at
the end of Section 9.5 on p.~118 of \cite{ca}, and easily proved
on lines 3--9 of \cite[p. 137]{ca}. Let us also mention that J\"oran
Bergh \cite{BerghJ1979} has obtained a stronger result related to
this inclusion.
\begin{rem}
The reader may care to compare Theorem \ref{babaeq} and its corollary
with Lemma 2 and its proof on pp.~66--68 of \cite{cj}.\textit{ }
\end{rem}
\textit{Proof of Theorem \ref{babaeq}.} First, for part (i), let
us suppose that $f\in\overline{\mathcal{F}}(X_{0}^{\#},X_{1}^{\#})$
with $\left\Vert f\right\Vert _{\overline{\mathcal{F}}(X_{0}^{\#},X_{1}^{\#})}=1$.
We consider the functions $f_{n}:\mathbb{S}\rightarrow X_{0}^{\#}+X_{1}^{\#}$
defined for each $n\in\mathbb{N}$ by 
\[
f_{n}(z)=\frac{n}{i}\left(f(z+i/n)-f(z)\right).
\]

Each function $f_{n}$ is of course a continuous map from $\mathbb{S}$
into $X_{0}^{\#}+X_{1}^{\#}$. We do not know a priori that $f_{n}$
is bounded. However, because of the definition of the space $\overline{\mathcal{F}}(X_{0}^{\#},X_{1}^{\#})$,
we know that 
\begin{equation}
\left\Vert f_{n}(z)\right\Vert _{X_{0}^{\#}+X_{1}^{\#}}\le c(1+\left|z\right|)\text{ for all }z\in\mathbb{S}\label{swx}
\end{equation}
for some constant $c$ which may depend on $n$. We also have $f_{n}(j+it)\in X_{j}^{\#}$
with 
\[
\left\Vert f_{n}(j+it)\right\Vert _{X_{j}^{\#}}\le1\text{ for }j=0,1\text{ and all }t\in\mathbb{R}.
\]

For each fixed $b\in X_{0}\cap X_{1}$ we consider the continuous
scalar function $\phi_{n,b}:\mathbb{S}\rightarrow\mathbb{C}$ defined
by $\phi_{n,b}(z)=\left\langle b,f_{n}(z)\right\rangle $. This is
clearly analytic in $\mathbb{S}^{\circ}$ and satisfies 
\begin{equation}
\left|\phi_{n,b}(j+it)\right|\le\left\Vert b\right\Vert _{X_{j}}\le\left\Vert b\right\Vert _{X_{0}\cap X_{1}}\text{ for }j=0,1\text{ and all }t\in\mathbb{R}.\label{hhl}
\end{equation}
 and 
\[
\left|\phi_{n,b}(z)\right|\le c\left\Vert b\right\Vert _{X_{0}\cap X_{1}}\left(1+\left|z\right|\right)\text{ for all }z\in\mathbb{S}.
\]
 So we can apply Lemma \ref{3l} to obtain that 
\begin{equation}
\left|\phi_{n,b}(z)\right|\le\left\Vert b\right\Vert _{X_{0}}^{1-\mathrm{Re\,}z}\left\Vert b\right\Vert _{X_{1}}^{\mathrm{Re\,}z}\le\left\Vert b\right\Vert _{X_{0}\cap X_{1}}\text{ for all }z\in\mathbb{S}^{\circ}.\label{prrqq}
\end{equation}
 Passing to the limit as $n$ tends to $\infty$ (and in fact using
(\ref{eq:wrnt})) we see that the $X_{0}^{\#}+X_{1}^{\#}$ valued
analytic function $h:=\left(\left.f\right|_{\mathbb{S}^{\circ}}\right)^{\prime}$
satisfies

\begin{equation}
\left|\left\langle b,h(z)\right\rangle \right|\le\left\Vert b\right\Vert _{X_{0}}^{1-\mathrm{Re\,}z}\left\Vert b\right\Vert _{X_{1}}^{\mathrm{Re\,}z}\le\left\Vert b\right\Vert _{X_{0}\cap X_{1}}\text{ for all }z\in\mathbb{S}^{\circ}.\label{rqq}
\end{equation}
So $h$ satisfies conditions (i) and (ii) of Definition \ref{upper}
with $C=1$. This establishes part (i) of Theorem \ref{babaeq}, except
that at this stage, instead of (\ref{wwx}), we only have 
\begin{equation}
\left\Vert f\right\Vert _{\overline{\mathcal{F}}(X_{0}^{\#},X_{1}^{\#})}\ge\left\Vert \left(\left.f\right|_{\mathbb{S}^{\circ}}\right)^{\prime}\right\Vert _{\mathcal{H}(X_{0}^{\#},X_{1}^{\#})}.\label{bawwx}
\end{equation}

Now, turning to the proof of part (ii) of the theorem, we suppose
that $h:\mathbb{S}\rightarrow X_{0}^{\#}+X_{1}^{\#}$ is an element
of $\mathcal{H}(X_{0}^{\#},X_{1}^{\#})$ with norm $1$. Then, for
each fixed $b\in X_{0}\cap X_{1}$, the function $h_{b}:\mathbb{S}^{\circ}\rightarrow\mathbb{C}$
defined by $h_{b}(z)=\left\langle b,h(z)\right\rangle $ is an element
of $H^{\infty}(\mathbb{S}^{\circ})$. It follows immediately from
condition (\ref{baconv}) that 
\begin{equation}
\left|h_{b}(z)\right|\le\left\Vert b\right\Vert _{X_{0}}^{1-\mathrm{Re}\, z}\left\Vert b\right\Vert _{X_{1}}^{\mathrm{Re}\, z}\le\left\Vert b\right\Vert _{X_{0}\cap X_{1}}^{1-\mathrm{Re}\, z}\left\Vert b\right\Vert _{X_{0}\cap X_{1}}^{\mathrm{Re}\, z}=\left\Vert b\right\Vert _{X_{0}\cap X_{1}}\text{ for all }z\in\mathbb{S}^{\circ}.\label{bral}
\end{equation}

Let $\phi_{b}:\mathbb{S}\rightarrow\mathbb{C}$ be the function obtained
from $h_{b}$ as in Lemma \ref{bapf}, by setting $\phi_{b}(z)=\int_{1/2}^{z}h_{b}(\zeta)d\zeta$
for all $z\in\mathbb{S}^{\circ}$ and then extending continuously
to $\mathbb{S}$. In view of condition (iv) of Lemma \ref{bapf} and
(\ref{bral}) we have 
\begin{equation}
\left|\phi_{b}(z)\right|\le\left|z-\frac{1}{2}\right|\cdot\left\Vert b\right\Vert _{X_{0}\cap X_{1}}\text{ for all }z\in\mathbb{S}.\label{bralz}
\end{equation}

\smallskip{}
Since $h_{b}(z)$ depends linearly on $b$ at each constant point
$z\in\mathbb{S}^{\circ}$, it follows that $\phi_{b}(z)$ also depends
linearly on $b$ at each constant point $z\in\mathbb{S}$. In view
of this and (\ref{bralz}), for each $z\in\mathbb{S}$ there exists
a unique element $f(z)\in(X_{0}\cap X_{1})^{\#}=X_{0}^{\#}+X_{1}^{\#}$
such that 
\begin{equation}
\left\langle b,f(z)\right\rangle =\phi_{b}(z)\text{ for all }b\in X_{0}\cap X_{1}\label{sdba}
\end{equation}
 and 
\begin{equation}
\left\Vert f(z)\right\Vert _{X_{0}^{\#}+X_{1}^{\#}}\le\left|z-\frac{1}{2}\right|.\label{ebta}
\end{equation}
 \smallskip{}
 Condition (i) of Lemma \ref{bapf} gives us that $\left|\phi_{b}(z_{1})-\phi_{b}(z_{2})\right|\le\left|z_{1}-z_{2}\right|\cdot\left\Vert b\right\Vert _{X_{0}\cap X_{1}}$
for all $b\in X_{0}\cap X_{1}$ and consequently $\left\Vert f(z_{1})-f(z_{2})\right\Vert _{X_{0}^{\#}+X_{1}^{\#}}\le\left|z_{1}-z_{2}\right|$
for all $z_{1},z_{2}\in\mathbb{S}$. This shows that $f:\mathbb{S}\rightarrow X_{0}^{\#}+X_{1}^{\#}$
is continuous.

From (\ref{sdba}) and the fact that $X_{0}\cap X_{1}$ is the predual
of $X_{0}^{\#}+X_{1}^{\#}$ and some standard results which are recalled
in Subsection \ref{hilphi} (cf.~also Remark \ref{rem:UseHilPhil})
we see that $f$ is an analytic $X_{0}^{\#}+X_{1}^{\#}$ valued function
on $\mathbb{S}^{\circ}$. Furthermore, it follows from (\ref{sdba})
that (cf.\ condition (ii) of Lemma \ref{bapf}) we must have $f^{\prime}(z)=h(z)$
for each $z\in\mathbb{S}^{\circ}$.

Condition (\ref{ebta}) and the discussion following it show that
$f$ satisfies the first three of the four conditions for membership
in the space $\overline{\mathcal{F}}(X_{0}^{\#},X_{1}^{\#})$ listed
in Section 5 on p.\ 115 of \cite{ca}. To obtain the fourth condition
we first deduce from (\ref{bral}) that (cf.\ Remark \ref{newntl})
$\lim_{r\searrow0}\left(\sup\left\{ \left|h_{b}(z)\right|:\left|\mathrm{Re\,}z-j\right|<r\right\} \right)\le\left\Vert b\right\Vert _{X_{j}}$
for $j=0,1$. Then condition (iii) of Lemma \ref{bapf} applied to
$h_{b}(z)$ tells us that 
\[
\left|\phi_{b}(j+it_{1})-\phi_{b}(j+it_{2})\right|\le\left|t_{1}-t_{2}\right|\cdot\left\Vert b\right\Vert _{X_{j}}\text{ }
\]
for $j=0,1$ and all $t_{1},t_{2}\in\mathbb{R}$ and all $b\in X_{0}\cap X_{1}$.

In view of this and (\ref{sdba}) it follows that 
\[
\left\Vert f\right\Vert _{\overline{\mathcal{F}}(X_{0}^{\#},X_{1}^{\#})}=\sup\left\{ \frac{\left\Vert f(j+it_{1})-f(j+it_{2})\right\Vert _{X_{j}^{\#}}}{|t_{1}-t_{2}|}:j=0,1,t_{1},t_{2}\in\mathbb{R},t_{1}\ne t_{2}\right\} \le1.
\]
This establishes the above-mentioned fourth condition, i.e., that
$f\in\overline{\mathcal{F}}(X_{0}^{\#},X_{1}^{\#})$, and it also
shows that in general 
\begin{equation}
\left\Vert f\right\Vert _{\overline{\mathcal{F}}(X_{0}^{\#},X_{1}^{\#})}\le\left\Vert h\right\Vert _{\mathcal{H}(X_{0}^{\#},X_{1}^{\#})}.\label{ilpe}
\end{equation}
But now, since $f^{\prime}=h$ on $\mathbb{S}^{\circ}$, we can also
apply the first part of the proof to obtain the reverse inequality
to (\ref{ilpe}), i.e., (\ref{bawwx}). This gives us (\ref{dlo}).
Analogously, in the context of the first part of this proof, we can
now obtain the reverse inequality to (\ref{bawwx}) for any given
$f\in\overline{\mathcal{F}}(X_{0}^{\#},X_{1}^{\#})$, by applying
the second part of the proof and (\ref{ilpe}) to the function $h\in\mathcal{H}(X_{0}^{\#},X_{1}^{\#})$
defined by $h=\left(\left.f\right|_{\mathbb{S}^{\circ}}\right)^{\prime}$.
This establishes (\ref{wwx}) and so completes the proof of Theorem
\ref{babaeq}. $\qed$

\smallskip{}

\subsection{\label{sub:OurProof}\protect Our proof of Theorem \ref{crt}.}

\hskip 1cm \phantom{A A}\\

\vskip -4mm \textit{Part 1: The inclusion} $[X_{0}^{\#},X_{1}^{\#}]^{\theta}\overset{1}{\subset}\left([X_{0},X_{1}]_{\theta}\right)^{\#}$\textit{.}

\smallskip{}

Suppose first that $y$ is an arbitrary element of $[X_{0}^{\#},X_{1}^{\#}]^{\theta}$
and that $\varepsilon$ is an arbitrary positive number. By Theorem
\ref{babaeq} and its corollary, there exists a continuous function
$h\in\mathcal{H}(X_{0}^{\#},X_{1}^{\#})$ such that $\left\Vert h\right\Vert _{\mathcal{H}(X_{0}^{\#},X_{1}^{\#})}\le(1+\varepsilon)\left\Vert y\right\Vert _{[X_{0}^{\#},X_{1}^{\#}]^{\theta}}$
and $h(\theta)=y$.

Now suppose that $x$ is an arbitrary element of $X_{0}\cap X_{1}$.
In view of Stafney's Lemma (\cite{stafney} Lemma 2.5 p.\ 335) there
exists a function $g$ in the space $\mathcal{G}_{\mathbf{St}}(X_{0},X_{1})$
(defined%
\footnote{We repeat our warning that Stafney's space $\mathcal{G}_{\mathbf{St}}(X_{0},X_{1})$
should not be confused with Calder\'on's smaller space $\mathcal{G}_{\mathbf{Ca}}\left(X_{0},X_{1}\right)$,
even though they are both denoted by $\mathcal{G}\left(X_{0},X_{1}\right)$
in \cite{stafney} and \cite{ca} respectively. Anyone who wishes
to apply Stafney's Lemma in other contexts should be aware that there
does not seem to be any obvious way of strengthening the proof of
that lemma to enable the function $g$ which satisfies $g(\theta)=x$
and (\ref{ggm}) to be assumed to be an element of $\mathcal{G}_{\mathbf{Ca}}\left(X_{0},X_{1}\right)$.
Other variants of Stafney's Lemma which can be used in conjunction
with various other interpolation methods can be found in \cite{Ivtsan}.%
} via Definition \ref{def:3DifferentGs} and \cite[p. 334]{stafney})
such that $g(\theta)=x$ and 
\begin{equation}
\sup_{t\in\mathbb{R}}\left\Vert g(j+it)\right\Vert _{X_{j}}\le\left(1+\frac{\varepsilon}{2}\right)\left\Vert x\right\Vert _{\left[X_{0},X_{1}\right]_{\theta}}\text{ for }j=0,1.\label{ggm}
\end{equation}

We recall that membership in $\mathcal{G}_{\mathbf{St}}(X_{0},X_{1})$
means that $g:\mathbb{S}\rightarrow X_{0}\cap X_{1}$ is a finite
sum of the form 
\begin{equation}
g(z)=\sum_{n=1}^{N}\phi_{n}(z)a_{n}\label{cbf}
\end{equation}
where each $a_{n}\in X_{0}\cap X_{1}$ and each $\phi_{n}:\mathbb{\mathbb{S}}\rightarrow\mathbb{C}$
is a continuous function on $\mathbb{S}$ which is analytic in $\mathbb{S}^{\circ}$
and satisfies $\lim_{R\rightarrow\infty}\sup\left\{ \left|\phi_{n}(z)\right|:z\in\mathbb{S},\left|\mathrm{Im\,}z\right|\ge R\right\} =0$.
Obviously 
\begin{equation}
\sup_{z\in\mathbb{S}}\left\Vert g(z)\right\Vert _{X_{j}}<\infty\mbox{ for }j=0,1,\label{eq:dfe}
\end{equation}
and furthermore, for $j=0,1$, the function $z\mapsto\left\Vert g(z)\right\Vert _{X_{j}}$
is uniformly continuous on $\mathbb{S}$. This latter property, together
with (\ref{ggm}) implies that, for some sufficiently small positive
number $r$, 
\begin{equation}
\sup\left\{ \left\Vert g(z)\right\Vert _{X_{j}}:z\in\mathbb{S},0<\left|j-\mathrm{Re\,\,}z\right|<r\right\} \le\left(1+\varepsilon\right)\left\Vert x\right\Vert _{\left[X_{0},X_{1}\right]_{\theta}}.\label{ppigv}
\end{equation}

\smallskip{}

We now consider the function $\Phi:\mathbb{S}^{\circ}\rightarrow\mathbb{C}$
defined by 
\[
\Phi(z)=\left\langle g(z),h(z)\right\rangle =\sum_{n=1}^{N}\phi_{n}(z)\left\langle a_{n},h(z)\right\rangle .
\]
This is a finite sum of products of functions in $H^{\infty}(\mathbb{S}^{\circ})$
and is therefore itself in $H^{\infty}(\mathbb{S}^{\circ})$. In view
of (\ref{baconv}) and (\ref{Aq}) in Definition \ref{upper}, we
also have

\begin{equation}
\left|\Phi(z)\right|\le(1+\varepsilon)\left\Vert y\right\Vert _{[X_{0}^{\#},X_{1}^{\#}]^{\theta}}\left\Vert g(z)\right\Vert _{X_{0}}^{1-\mathrm{Re\,}z}\left\Vert g(z)\right\Vert _{X_{1}}^{\mathrm{Re\,}z}\text{ for all }z\in\mathbb{S}^{\circ}.\label{isms}
\end{equation}

Using (\ref{isms}), (\ref{ppigv}) and (\ref{eq:dfe}), we obtain
that 
\[
\lim_{r\searrow0}\,\sup\left\{ \left|\Phi(z)\right|:z\in\mathbb{S},\text{ }0<\left|j-\mathrm{Re\,}z\right|<r\right\} \le(1+\varepsilon)^{2}\left\Vert x\right\Vert _{\left[X_{0},X_{1}\right]_{\theta}}\left\Vert y\right\Vert _{[X_{0}^{\#},X_{1}^{\#}]^{\theta}}\text{ for }j=0,1.
\]
 This implies, using also (\ref{eq:dfe}) again, together with part
(ii) of Lemma \ref{3l}, that 
\[
\left|\left\langle x,y\right\rangle \right|=\left|\left\langle g(\theta),h(\theta)\right\rangle \right|=\left|\Phi(\theta)\right|\le(1+\varepsilon)^{2}\left\Vert x\right\Vert _{\left[X_{0},X_{1}\right]_{\theta}}\left\Vert y\right\Vert _{[X_{0}^{\#},X_{1}^{\#}]^{\theta}}.
\]
Since $\varepsilon$ can be chosen arbitrarily small, we have shown
that
\begin{equation}
\left|\left\langle x,y\right\rangle \right|\le\left\Vert x\right\Vert _{\left[X_{0},X_{1}\right]_{\theta}}\left\Vert y\right\Vert _{[X_{0}^{\#},X_{1}^{\#}]^{\theta}}\mbox{ for each }x\in X_{0}\cap X_{1}\mbox{ and each }y\in\left[X_{0}^{\#},X_{1}^{\#}\right]^{\theta}.\label{eq:PreStuff}
\end{equation}
Consequently (cf.~Definition \ref{blg}), each $y\in\left[X_{0}^{\#},X_{1}^{\#}\right]^{\theta}$
is also an element of $\left([X_{0},X_{1}]_{\theta}\right)^{\#}$
and satisfies $\left\Vert y\right\Vert _{\left([X_{0},X_{1}]_{\theta}\right)^{\#}}\le\left\Vert y\right\Vert _{[X_{0}^{\#},X_{1}^{\#}]^{\theta}}$. 

In other words, we have shown that $[X_{0}^{\#},X_{1}^{\#}]^{\theta}\overset{1}{\subset}\left([X_{0},X_{1}]_{\theta}\right)^{\#}$
and therefore completed the proof of Part 1. 

But let us also remark that (\ref{eq:PreStuff}) and the fact that
$X_{0}\cap X_{1}$ is dense in $\left[X_{0},X_{1}\right]_{\theta}$
immediately imply (cf.~Definition \ref{ranX}) that 
\begin{equation}
\left|\left\langle x,y\right\rangle _{\left[X_{0},X_{1}\right]_{\theta}}\right|\le\left\Vert x\right\Vert _{\left[X_{0},X_{1}\right]_{\theta}}\left\Vert y\right\Vert _{[X_{0}^{\#},X_{1}^{\#}]^{\theta}}\mbox{ for each }x\in\left[X_{0},X_{1}\right]_{\theta}\mbox{ and each }y\in\left[X_{0}^{\#},X_{1}^{\#}\right]^{\theta}.\label{eq:fatu}
\end{equation}

\smallskip{}

\textit{Part 2: The inclusion} $\left([X_{0},X_{1}]_{\theta}\right)^{\#}\overset{1}{\subset}[X_{0}^{\#},X_{1}^{\#}]^{\theta}$\textit{.}

\smallskip{}

\smallskip{}
Let us now suppose that $y$ is an arbitrary element of $\left([X_{0},X_{1}]_{\theta}\right)^{\#}$.
We shall begin by using $y$ to define a linear functional $\lambda$
on the space $\mathcal{G}_{\mathbb{D}}(X_{0},X_{1})$ (see Definition
\ref{badefg}) by setting $\lambda(f)=\left\langle f(0),y\right\rangle $
for each $f\in\mathcal{G}_{\mathbb{D}}(X_{0},X_{1})$.

(Note that $f(0)\in X_{0}\cap X_{1}$ whenever $f\in\mathcal{G}_{\mathbb{D}}(X_{0},X_{1})$,
and, (cf.\ Fact \ref{fkz} and (\ref{zaq})) we also have $\left([X_{0},X_{1}]_{\theta}\right)^{\#}\subset\left(X_{0}\cap X_{1}\right)^{\#}$
so the notation $\left\langle f(0),y\right\rangle $ is appropriate.)
Instead of using the norm (\ref{tulfamw}), we shall equip $\mathcal{G}_{\mathbb{D}}(X_{0},X_{1})$
with the (smaller) norm 
\begin{equation}
\left\Vert f\right\Vert _{L}=\frac{1}{2\pi}\left(\int_{2\pi\theta}^{2\pi}\left\Vert f(e^{it})\right\Vert _{X_{0}}dt+\int_{0}^{2\pi\theta}\left\Vert f(e^{it})\right\Vert _{X_{1}}dt\right)\label{uypl}
\end{equation}
 In view of (\ref{aqz}), we see that $\lambda$ is bounded on $\left(\mathcal{G}_{\mathbb{D}}(X_{0},X_{1}),\left\Vert \cdot\right\Vert _{L}\right)$
with norm not exceeding $\left\Vert y\right\Vert _{\left([X_{0},X_{1}]_{\theta}\right)^{\#}}$. 

Let $L$ be the space of (equivalence classes of) functions $f:\mathbb{T}\rightarrow X_{0}\cap X_{1}$
which are finite sums of the form $f(w)=\sum_{k=1}^{N}\psi_{k}(w)a_{k}$,
where, again, as in Definition \ref{badefg}, each $a_{k}$ is in
$X_{0}\cap X_{1}$, but now each $\psi_{k}:\mathbb{T}\rightarrow\mathbb{C}$
is only defined on $\mathbb{T}$, and $\psi(e^{it})$ is an integrable
function of $t$ on $[0,2\pi]$. We use the same norm (\ref{uypl})
for $L$ as we did for $\mathcal{G}_{\mathbb{D}}(X_{0},X_{1})$. It
is easy to see that, for each $f\in L$, the functions $t\mapsto\left\Vert f(e^{it})\right\Vert _{X_{0}}$
and $t\mapsto\left\Vert f(e^{it})\right\Vert _{X_{1}}$ are both Lebesgue
measurable. (One can deduce this from the fact that $f$ takes values
in a finite dimensional subspace of $X_{0}\cap X_{1}$. Alternatively,
one can observe that each scalar function $\psi_{k}$ is the pointwise
limit of a sequence of simple functions, and that therefore the same
is true for $t\mapsto\left\Vert f(e^{it})\right\Vert _{X_{j}}$.)
Therefore the two integrals in (\ref{uypl}) are well defined.

Since, for each $f\in\mathcal{G}_{\mathbb{D}}(X_{0},X_{1})$ the values
taken by $f$ on $\mathbb{T}$ determine its values on all of $\mathbb{D}$
and also its norm $\left\Vert f\right\Vert _{L}$, we can consider
$\left(\mathcal{G}_{\mathbb{D}}(X_{0},X_{1}),\left\Vert \cdot\right\Vert _{L}\right)$
to be a space of functions on $\mathbb{T}$ and thus it is a subspace
of $L$. (In more precise language, we are dealing here with the space
which consists of every function $f:\mathbb{T}\to X_{0}\cap X_{1}$
which has a (necessarily unique) extension to $\mathbb{D}$ which
is an element of $\mathcal{G}_{\mathbb{D}}(X_{0},X_{1})$.) The Hahn--Banach
theorem guarantees the existence of a bounded linear functional $\widetilde{\lambda}:L\rightarrow\mathbb{C}$
which is an extension of $\lambda$ and which satisfies $\left|\widetilde{\lambda}(f)\right|\le\left\Vert y\right\Vert _{\left([X_{0},X_{1}]_{\theta}\right)^{\#}}\cdot\left\Vert f\right\Vert _{L}$
for each $f\in L$. Let use the usual notation $L^{1}(\mathbb{T})$
to denote the familiar space of all (equivalence classes of) measurable
functions $\phi:\mathbb{T}\rightarrow\mathbb{C}$ for which $\left\Vert \phi\right\Vert _{L^{1}(\mathbb{T})}:=\frac{1}{2\pi}\int_{0}^{2\pi}\left|\phi(e^{it})\right|dt<\infty$.
For each fixed $b\in X_{0}\cap X_{1}$, we define a linear functional
$\widetilde{\lambda}_{b}$ on this space by the formula 
\begin{equation}
\widetilde{\lambda}_{b}(\phi)=\widetilde{\lambda}(\phi b)\label{pfoza}
\end{equation}
for all $\phi\in L^{1}(\mathbb{T})$. Here of course $\phi b$ means
the function in $L$ defined by $(\phi b)(e^{it})=\phi(e^{it})b$.
Thus, for each $\phi\in L^{1}(\mathbb{T})$, we have 
\begin{eqnarray*}
\left|\widetilde{\lambda}_{b}(\phi)\right| & \le & \left\Vert y\right\Vert _{\left([X_{0},X_{1}]_{\theta}\right)^{\#}}\cdot\left\Vert \phi b\right\Vert _{L}\\
 & = & \left\Vert y\right\Vert _{\left([X_{0},X_{1}]_{\theta}\right)^{\#}}\cdot\frac{1}{2\pi}\left(\left\Vert b\right\Vert _{X_{0}}\int_{2\pi\theta}^{2\pi}\left|\phi(e^{it})\right|dt+\left\Vert b\right\Vert _{X_{1}}\int_{0}^{2\pi\theta}\left|\phi(e^{it})\right|dt\right)\\
 & = & \left\Vert y\right\Vert _{\left([X_{0},X_{1}]_{\theta}\right)^{\#}}\cdot\frac{1}{2\pi}\int_{0}^{2\pi}w(t)\left|\phi(e^{it})\right|dt,
\end{eqnarray*}
where $w=\left\Vert b\right\Vert _{X_{0}}\chi_{[2\pi\theta,2\pi]}+\left\Vert b\right\Vert _{X_{1}}\chi_{[0,2\pi\theta]}$.

It is easy to deduce, from the standard isometric identification of
$\left(L^{1}(\mathbb{T})\right)^{*}$ with $L^{\infty}(\mathbb{T})$,
that there exists an essentially bounded measurable function $h_{b}:\mathbb{T}\rightarrow\mathbb{C}$
such that 
\begin{equation}
\widetilde{\lambda}_{b}(\phi)=\frac{1}{2\pi}\int_{0}^{2\pi}\phi(e^{it})h_{b}(e^{it})dt\text{ for all }\phi\in L^{1}(\mathbb{T})\text{ }\label{foza}
\end{equation}
and $\underset{t\in[0,2\pi]}{\mathrm{ess~sup}}\left|h_{b}(e^{it})/w(t)\right|\le\left\Vert y\right\Vert _{\left([X_{0},X_{1}]_{\theta}\right)^{\#}}$,
which we can also write as 
\begin{equation}
\underset{e^{it}\in\Gamma_{j}}{\mathrm{ess~sup}}\left|h_{b}(e^{it})\right|\le\left\Vert y\right\Vert _{\left([X_{0},X_{1}]_{\theta}\right)^{\#}}\cdot\left\Vert b\right\Vert _{X_{j}}\text{ for }j=0,1.\label{mskgw}
\end{equation}

\smallskip{}

For each positive constant $\delta>0$, the function $\phi:\mathbb{S}\to\mathbb{C}$
defined by $\phi(z)=e^{\delta z^{2}}$ satisfies $\sup_{z\in\mathbb{S}}\left|\phi(z)\right|\le e^{\delta}$.
Since it also satisfies the conditions of part (i) of Lemma \ref{lem:PhiPsi}
we see that the function $\Omega_{\delta}:\mathbb{D}\rightarrow\mathbb{C}$
defined by 
\[
\Omega_{\delta}(w):=\left\{ \begin{array}{lll}
e^{\delta\left(\xi(w)\right)^{2}} & , & w\in\mathbb{D}\backslash\{1,e^{i2\pi\theta}\}\\
0 & , & w=1,e^{i2\pi\theta}
\end{array}\right.
\]
is continuous on $\mathbb{D}$ and analytic in $\mathbb{D}^{\circ}$.
Consequently, for each $k\in\mathbb{N}$ and each $b\in X_{0}\cap X_{1}$,
the function $f_{\delta,k,b}(w):=w^{k}\Omega_{\delta}(w)b$ is an
element of $\mathcal{G}_{\mathbb{D}}(X_{0},X_{1})$. So we have 
\begin{eqnarray*}
0 & = & \left\langle f_{\delta,k,b}(0),y\right\rangle =\lambda(f_{\delta,k,b})=\widetilde{\lambda}(f_{\delta,k,b})=\widetilde{\lambda}_{b}(w^{k}\Omega_{\delta}(w))\\
 & = & \frac{1}{2\pi}\int_{0}^{2\pi}e^{ikt}e^{\delta\left(\xi(e^{it})\right)^{2}}h_{b}(e^{it})dt.
\end{eqnarray*}
The fact that $\sup_{w\in\mathbb{D}}\left|\Omega_{\delta}(w)\right|=\sup_{z\in\mathbb{S}}\left|e^{\delta z^{2}}\right|\le e^{\delta}$
enables us to apply the dominated convergence theorem in order to
take the limit as $\delta$ tends to $0$ of the preceding integral.
This gives us that $\frac{1}{2\pi}\int_{0}^{2\pi}e^{ikt}h_{b}(e^{it})dt=0$
for all $k\in\mathbb{N}$. I.e., the function $h_{b}$ satisfies the
condition in part (iii) of Theorem \ref{pfatou}, and obviously also
(\ref{mybigp}). Therefore the function $u_{b}:\mathbb{D}^{\circ}\rightarrow\mathbb{C}$
defined for each $w\in\mathbb{D}^{\circ}$ by the absolutely convergent
series (cf.\ (\ref{defu})) 
\[
u_{b}(w):=\sum_{n=1}^{\infty}\left(\frac{1}{2\pi}\int_{0}^{2\pi}e^{-ins}h_{b}(e^{is})ds\right)w^{n}
\]
is analytic in $\mathbb{D}^{\circ}$. Part (iv) of Theorem \ref{pfatou},
together with (\ref{mskgw}) gives us that 
\begin{equation}
\left|u_{b}(w)\right|\le\left\Vert y\right\Vert _{\left([X_{0},X_{1}]_{\theta}\right)^{\#}}\cdot\left\Vert b\right\Vert _{X_{0}\cap X_{1}}\text{ for each }w\in\mathbb{D}^{\circ}\text{.}\label{wengw}
\end{equation}
Part (vi) of Theorem \ref{pfatou}, together with (\ref{mskgw}) gives
us that the bounded analytic function $\phi_{b}:\mathbb{S}^{\circ}\rightarrow\mathbb{C}$
defined by $\phi_{b}(z)=u_{b}(\xi^{-1}(z))$ satisfies 
\begin{equation}
\left|\phi_{b}(z)\right|=\left|u_{b}(\xi^{-1}(z))\right|\le\left\Vert y\right\Vert _{\left([X_{0},X_{1}]_{\theta}\right)^{\#}}\cdot\left\Vert b\right\Vert _{X_{0}}^{1-\mathrm{Re\,}z}\cdot\left\Vert b\right\Vert _{X_{1}}^{\mathrm{Re\,}z}\text{ for all }z\in\mathbb{S}^{\circ}.\label{wanegw}
\end{equation}
Each Fourier coefficient $\frac{1}{2\pi}\int_{0}^{2\pi}e^{-ins}h_{b}(e^{is})ds$
of $h_{b}$ depends linearly on $b$ since, when we apply (\ref{foza})
and then (\ref{pfoza}) to the function $\phi\in L^{1}(\mathbb{T})$
defined by $\phi(e^{it})=e^{-int}$, we have 
\[
\frac{1}{2\pi}\int_{0}^{2\pi}e^{-ins}h_{b}(e^{is})ds=\widetilde{\lambda}_{b}(\phi)=\widetilde{\lambda}(\phi b).
\]
It follows that $u_{b}(w)$ and $\phi_{b}(z)$ also depend linearly
on $b$, for each fixed $w\in\mathbb{D}^{\circ}$ and each fixed $z\in\mathbb{S}^{\circ}$.
Thus, by (\ref{wengw}) these are both bounded linear functionals
on $X_{0}\cap X_{1}$, and, for each $z\in\mathbb{S}^{\circ}$, there
exists an element $v(z)\in(X_{0}\cap X_{1})^{\#}=X_{0}^{\#}+X_{1}^{\#}$
such that $\left\langle b,v(z)\right\rangle =\phi_{b}(z)$ for each
$b\in X_{0}\cap X_{1}$. The properties of $\phi_{b}$ stated above,
notably (\ref{wanegw}), imply that $v\in\mathcal{H}(X_{0}^{\#},X_{1}^{\#})$
with $\left\Vert v\right\Vert _{\mathcal{H}(X_{0}^{\#},X_{1}^{\#})}\le\left\Vert y\right\Vert _{\left([X_{0},X_{1}]_{\theta}\right)^{\#}}$.
Thus it remains only to show that $v(\theta)=y$, i.e., that $\left\langle b,v(\theta)\right\rangle =\left\langle b,y\right\rangle $
for each $b\in X_{0}\cap X_{1}$. In view of the various definitions
given in the preceding steps, we see that

\begin{eqnarray*}
\left\langle b,v(\theta)\right\rangle  & = & \phi_{b}(\theta)=u_{b}(\xi^{-1}(\theta))=u_{b}(0)=\frac{1}{2\pi}\int_{0}^{2\pi}h_{b}(e^{is})ds\\
 & = & \lim_{\delta\searrow0}\frac{1}{2\pi}\int_{0}^{2\pi}e^{\delta\xi(e^{it})^{2}}h_{b}(e^{is})ds\\
 & = & \lim_{\delta\searrow0}\widetilde{\lambda}_{b}\left(\Omega_{\delta}\right)=\lim_{\delta\searrow0}\widetilde{\lambda}\left(\Omega_{\delta}b\right)=\lim_{\delta\searrow0}\lambda\left(\Omega_{\delta}b\right)\\
 & = & \lim_{\delta\searrow0}\left\langle \Omega_{\delta}(0)b,y\right\rangle =\left\langle b,y\right\rangle \text{ ,}
\end{eqnarray*}
 which indeed completes the proof. $\qed$ 

\medskip{}

\begin{rem}
\label{rem:InfAttained}We can use the steps of the preceding proof
to show that the infimum (\ref{eq:wfnt}) in the alternative definition
in Corollary \ref{corbab} of the norm in the space $\left[X_{0}^{\#},X_{1}^{\#}\right]^{\theta}$
is always attained, for every Banach couple $\left(X_{0}^{\#},X_{1}^{\#}\right)$
which arises as the dual couple of a regular couple $\left(X_{0},X_{1}\right)$.
In fact, via Theorem \ref{babaeq}, or via an examination of the original
proof of Theorem \ref{crt} in \cite{ca} or \cite{bl}, the analogous
infimum in the original definition of $\left[X_{0}^{\#},X_{1}^{\#}\right]^{\theta}$
is also always attained for such couples. More explicitly, given any
element $y\in\left[X_{0}^{\#},X_{1}^{\#}\right]^{\theta}$ we can
first use Part 1 of the above proof to show that $y$ is also an element
of the dual space $\left(\left[X_{0},X_{1}\right]_{\theta}\right)^{*}$
and satisfies $\left\Vert y\right\Vert _{\left([X_{0},X_{1}]_{\theta}\right)^{\#}}\le\left\Vert y\right\Vert _{[X_{0}^{\#},X_{1}^{\#}]^{\theta}}$.
But then, considering $y$ as an element of the dual space, the steps
of Part 2 of the same proof show that there exists a function $v\in\mathcal{H}(X_{0}^{\#},X_{1}^{\#})$
with $\left\Vert v\right\Vert _{\mathcal{H}(X_{0}^{\#},X_{1}^{\#})}\le\left\Vert y\right\Vert _{\left([X_{0},X_{1}]_{\theta}\right)^{\#}}$
and $v(\theta)=y$. By definition, we must then also have $\left\Vert y\right\Vert _{\left[X_{0}^{\#},X_{1}^{\#}\right]^{\theta}}\le\left\Vert v\right\Vert _{\mathcal{H}(X_{0}^{\#},X_{1}^{\#})}$.
So all inequalities here have to be equalities and $v$ attains the
infimum in (\ref{eq:wfnt}).
\end{rem}

\section{Some consequences of Calder\'on's duality theorem}

\subsection{\label{sub:UsefulNormingSubspace}A useful norming subspace of $\left[X_{0},X_{1}\right]_{\theta}^{\#}$}

Obviously, having a concrete description of the dual of $\left[X_{0},X_{1}\right]_{\theta}$
provides us with a powerful tool for the study of this space. But
for some purposes it turns out to be more convenient not to work with
the whole of the dual space $\left[X_{0}^{\#},X_{1}^{\#}\right]^{\theta}$.
The definition of $\left[X_{0}^{\#},X_{1}^{\#}\right]^{\theta}$ is
relatively complicated. Fortunately some things can be done by using
its apparently simpler subspace $\left[X_{0}^{\#},X_{1}^{\#}\right]_{\theta}$
or even the dense subspace $X_{0}^{\#}\cap X_{1}^{\#}$ of $\left[X_{0}^{\#},X_{1}^{\#}\right]_{\theta}$.
For example, in \cite{ciddar}, properties of these subspaces enable
Calderón's reiteration formula for his interpolation spaces to be
proved without need for the auxiliary condition which was imposed
in \cite{ca}. (An alternative and quite short proof of this formula,
but only to within equivalence of norms, is given by Janson in \cite{janson}.)
Then, in \cite{cwsh}, these same subspaces help to give an isometric
version of a result of Janson \cite{janson} about complex interpolation
of Gagliardo closures of given Banach spaces. 

Here is the theorem (cf.~\cite{ciddar}) which makes it possible
to use $\left[X_{0}^{\#},X_{1}^{\#}\right]_{\theta}$ and $X_{0}^{\#}\cap X_{1}^{\#}$
in such contexts. In its formulation we use the bilinear functional
$\left\langle \cdot,\cdot\right\rangle _{X}$ introduced in Definition
\ref{ranX}. Although we have formulated it for regular Banach couples,
it is of course also immediately applicable to arbitrary Banach couples
by passing to their regularizations (cf.~Remark \ref{rem:IfNotRegular}). 
\begin{thm}
\label{thm:oldbaba}Let $\left(X_{0},X_{1}\right)$ be an arbitrary
regular Banach couple. For each $\theta\in(0,1)$, for each $a\in\left[X_{0},X_{1}\right]_{\theta}$
and for each $\varepsilon>0$, there exists an element $b\in X_{0}^{\#}\cap X_{1}^{\#}$
such that
\begin{equation}
(1-\varepsilon)\left\Vert a\right\Vert _{\left[X_{0},X_{1}\right]_{\theta}}\cdot\left\Vert b\right\Vert _{\left[X_{0}^{\#},X_{1}^{\#}\right]_{\theta}}\le\left|\left\langle a,b\right\rangle _{X_{0}+X_{1}}\right|=\left|\left\langle a,b\right\rangle _{[X_{0},X_{1}]_{\theta}}\right|\le\left\Vert a\right\Vert _{\left[X_{0},X_{1}\right]_{\theta}}\cdot\left\Vert b\right\Vert _{\left[X_{0}^{\#},X_{1}^{\#}\right]_{\theta}}.\label{eq:faab}
\end{equation}

\end{thm}
We will present the proof of this theorem later, immediately after
Lemma \ref{lem:Preparation} in Subsection \ref{sub:AnotherLook},
since that lemma supplies the main ingredient that we need for the
proof.

\subsection{\label{sub:AnotherLook}Another look at page 1006 of \cite{ciddar}}

Let us here, as promised above, rewrite some of the arguments of the
lower half of page 1006 of \cite{ciddar} more pedantically, being
explicit about which duality (or bilinear functional) we are using
at each stage. This requires us to use notation and notions from Section
\ref{nbd}. We will also need some notation and notions from Section
\ref{cald}. Along the way we will also collect a result which will
enable us to supply the proof of Theorem \ref{thm:oldbaba}. 

Let $\left(A_{0},A_{1}\right)$ be an arbitrary Banach couple of complex
Banach spaces and choose $\theta\in\left(0,1\right)$. Since $\left(A_{0},A_{1}\right)$
is not necessarily regular we let $\left(X_{0},X_{1}\right)$ be the
regularization of $\left(A_{0},A_{1}\right)$ (cf.~Remark \ref{rem:IfNotRegular})
and keep in mind that $\left[A_{0},A_{1}\right]_{\theta}\overset{1}{=}\left[X_{0},X_{1}\right]_{\theta}$.

The following lemma expresses the particular result on page 1006 of
\cite{ciddar}, whose proof we wish to rewrite here.
\begin{lem}
\label{lem:3First}For each bounded linear functional $\ell\in\left(\left[A_{0},A_{1}\right]_{\theta}\right)^{*}$,
the action of $\ell$ on $\left[A_{0},A_{1}\right]_{\theta}$ can
be expressed by the formula 
\end{lem}
\begin{equation}
\ell(a)=\lim_{r\searrow0}\left\langle a,h_{r}(\theta)\right\rangle _{[X_{0},X_{1}]_{\theta}}\mbox{ for all }a\in[X_{0},X_{1}]_{\theta}\label{eq:oaists}
\end{equation}
for a certain family of functions $\left\{ h_{r}\right\} _{r>0}$
in $\mathcal{F}(X_{0}^{\#},X_{1}^{\#})$ depending on $\ell$, which
satisfy 
\begin{equation}
\left\Vert h_{r}\right\Vert _{\mathcal{F}(X_{0}^{\#},X_{1}^{\#})}\le e^{r}\left\Vert \ell\right\Vert _{\left(\left[X_{0},X_{1}\right]_{\theta}\right)^{*}}\mbox{ for each }r>0\label{eq:afyd}
\end{equation}
(Here we have slightly changed the notation of \cite{ciddar} by writing
$\theta$ instead of $s$, and $r$ instead of $\eta$.)

\noindent\textit{Proof of Lemma \ref{lem:3First}.} Let $\ell$ be
an arbitrary bounded linear functional acting on $\left[A_{0},A_{1}\right]_{\theta}$.
The Calder\'on duality theorem, i.e., Theorem \ref{crt}, (no matter
how we prefer to prove it), augmented by Remark \ref{rem:InfAttained},
tells us that there exists a function $h\in\overline{\mathcal{F}}(X_{0}^{\#},X_{1}^{\#})$
with $\left\Vert h\right\Vert _{\overline{\mathcal{F}}(X_{0}^{\#},X_{1}^{\#})}=\left\Vert \ell\right\Vert _{\left(\left[X_{0},X_{1}\right]_{\theta}\right)^{*}}$
and $\left\Vert h'(\theta)\right\Vert _{\left[X_{0}^{\#},X_{1}^{\#}\right]^{\theta}}=\left\Vert \ell\right\Vert _{\left(\left[X_{0},X_{1}\right]_{\theta}\right)^{*}}$
for which 
\[
\ell(a)=\left\langle a,h'(\theta)\right\rangle =\left\langle a,h'(\theta)\right\rangle _{X_{0}\cap X_{1}}=\left\langle a,h'(\theta)\right\rangle _{[X_{0},X_{1}]_{\theta}}\mbox{ for each }a\in A_{0}\cap A_{1}=X_{0}\cap X_{1}.
\]
(The last of these equalities is simply an application of Definition
\ref{ranX} for $X=\left[X_{0},X_{1}\right]_{\theta}$.) We can now
specify that the functions $h_{r}$ which appear in (\ref{eq:oaists})
are to be obtained from this function $h$, by setting $h_{r}(z):=e^{rz^{2}}\left(h(z+ir\right)-h(z))/ir$
for all $z\in\mathbb{S}$ and $r>0$. So, (in view of (\ref{eq:wrnt})
and, more specifically, (\ref{eq:vvdia}) of Remark \ref{rem:WeakerAnalyticity})
we have that 
\begin{equation}
\left\langle a,h'(\theta)\right\rangle _{X_{0}\cap X_{1}}=\lim_{r\searrow0}\left\langle a,h_{r}(\theta)\right\rangle _{X_{0}\cap X_{1}}\mbox{ for each }a\in X_{0}\cap X_{1}.\label{eq:emdalia}
\end{equation}
(We can observe that this proves (\ref{eq:oaists}) for the special
case where $a\in A_{0}\cap A_{1}$. But we need to deal with the general
case.) 

It is obvious (as already observed in \cite[p. 136]{ca} and in \cite{ciddar})
that, for each $r>0$, we have that $h_{r}\in\mathcal{F}(X_{0}^{\#},X_{1}^{\#})$
with $\left\Vert h_{r}\right\Vert _{\mathcal{F}(X_{0}^{\#},X_{1}^{\#})}\le e^{r}\left\Vert \ell\right\Vert _{\left(\left[X_{0},X_{1}\right]_{\theta}\right)^{*}}$,
for each $r>0$. Therefore we also have $h_{r}(\theta)\in[X_{0}^{\#},X_{1}^{\#}]_{\theta}$
with $\left\Vert h_{r}(\theta)\right\Vert _{\left[X_{0}^{\#},X_{1}^{\#}\right]_{\theta}}\le e^{r}\left\Vert \ell\right\Vert _{\left(\left[X_{0},X_{1}\right]_{\theta}\right)^{*}}$. 

The inequality (\ref{eq:fatu}) obtained in the course of proving
Theorem \ref{crt}, when combined with the norm one inclusion (\ref{eq:Norm1Inclusion})
immediately implies that 
\begin{equation}
\left|\left\langle x,y\right\rangle _{\left[X_{0},X_{1}\right]_{\theta}}\right|\le\left\Vert x\right\Vert _{\left[X_{0},X_{1}\right]_{\theta}}\left\Vert y\right\Vert _{\left[X_{0}^{\#},X_{1}^{\#}\right]_{\theta}}\mbox{ for each }x\in\left[X_{0},X_{1}\right]_{\theta}\mbox{ and }y\in\left[X_{0}^{\#},X_{1}^{\#}\right]_{\theta}.\label{eq:frtu}
\end{equation}
(We should perhaps recall that the definition of $\left\langle x,y\right\rangle _{\left[X_{0},X_{1}\right]_{\theta}}$
is as specified by Definition \ref{ranX}.)

Now let $a$ be an arbitrary element of $\left[X_{0},X_{1}\right]_{\theta}$.
Given $\varepsilon>0$, there exists an element $b_{\varepsilon}\in X_{0}\cap X_{1}$
such that 
\begin{equation}
\left\Vert a-b_{\varepsilon}\right\Vert _{\left[X_{0},X_{1}\right]_{\theta}}\le\varepsilon.\label{eq:fvtu}
\end{equation}
Then, for each $r>0$, we have 
\begin{eqnarray*}
 &  & \left|\ell(a)-\left\langle a,h_{r}\left(\theta\right)\right\rangle _{\left[X_{0},X_{1}\right]_{\theta}}\right|\\
 & = & \left|\left\langle a,h'\left(\theta\right)\right\rangle _{\left[X_{0},X_{1}\right]_{\theta}}-\left\langle a,h_{r}\left(\theta\right)\right\rangle _{\left[X_{0},X_{1}\right]_{\theta}}\right|\\
 & = & \left|\left\langle a-b_{\varepsilon},h^{\prime}(\theta)\right\rangle _{\left[X_{0},X_{1}\right]_{\theta}}+\left\langle b_{\varepsilon},h^{\prime}(\theta)-h_{r}(\theta)\right\rangle _{\left[X_{0},X_{1}\right]_{\theta}}+\left\langle b_{\varepsilon}-a,h_{r}(\theta)\right\rangle _{\left[X_{0},X_{1}\right]_{\theta}}\right|\\
 & \le & \left|\left\langle a-b_{\varepsilon},h^{\prime}(\theta)\right\rangle _{\left[X_{0},X_{1}\right]_{\theta}}\right|+\left|\left\langle b_{\varepsilon},h^{\prime}(\theta)-h_{r}(\theta)\right\rangle _{\left[X_{0},X_{1}\right]_{\theta}}\right|+\left|\left\langle b_{\varepsilon}-a,h_{r}(\theta)\right\rangle _{\left[X_{0},X_{1}\right]_{\theta}}\right|.
\end{eqnarray*}
In view of (\ref{eq:fatu}), (\ref{eq:frtu}) and (\ref{eq:fvtu})
and the above estimates for $\left\Vert h\right\Vert _{\overline{\mathcal{F}}(X_{0}^{\#},X_{1}^{\#})}$
and $\left\Vert h_{r}\right\Vert _{\mathcal{F}(X_{0}^{\#},X_{1}^{\#})}$,
this last expression does not exceed

\begin{eqnarray*}
 &  & \left\Vert a-b_{\varepsilon}\right\Vert _{\left[X_{0},X_{1}\right]_{\theta}}\left\Vert h'\left(\theta\right)\right\Vert _{\left[X_{0}^{\#},X_{1}^{\#}\right]^{\theta}}+\left|\left\langle b_{\varepsilon},h^{\prime}(\theta)-h_{r}(\theta)\right\rangle _{\left[X_{0},X_{1}\right]_{\theta}}\right|+\left\Vert a-b_{\varepsilon}\right\Vert _{\left[X_{0},X_{1}\right]_{\theta}}e^{r}\left\Vert \ell\right\Vert _{\left(\left[X_{0},X_{1}\right]_{\theta}\right)^{*}}\\
 & \le & \varepsilon\left\Vert \ell\right\Vert _{\left(\left[X_{0},X_{1}\right]_{\theta}\right)^{*}}+\left|\left\langle b_{\varepsilon},h^{\prime}(\theta)-h_{r}(\theta)\right\rangle _{\left[X_{0},X_{1}\right]_{\theta}}\right|+\varepsilon e^{r}\left\Vert \ell\right\Vert _{\left(\left[X_{0},X_{1}\right]_{\theta}\right)^{*}}.
\end{eqnarray*}
In view of (\ref{eq:emdalia}) the middle term tends to $0$ as $r\to0$.
Thus we have shown that 
$$\lim_{r\searrow0}
\left|\ell(a)-\left\langle a,h_{r}\left(\theta\right)\right\rangle _{\left[X_{0},X_{1}\right]_{\theta}}\right|
\le2\varepsilon\left\Vert \ell\right\Vert _{\left(\left[X_{0},X_{1}\right]_{\theta}\right)^{*}}$$
for our particular (and arbitrary) choice of positive $\varepsilon$.
This completes the proof of (\ref{eq:oaists}) and so also of Lemma
\ref{lem:3First}. $\qed$ 
\begin{rem}
This is a convenient place to mention and prove the following simple
consequence of (\ref{eq:oaists}) which will then enable us to promptly
prove Theorem \ref{thm:oldbaba} of Subsection \ref{sub:UsefulNormingSubspace}.
\end{rem}

\begin{lem}
\label{lem:Preparation} Let $\left(A_{0},A_{1}\right)$ be an arbitrary
Banach couple and let $\left(X_{0},X_{1}\right)$ be the regularization
of $\left(A_{0},A_{1}\right)$ (cf.~Remark \ref{rem:IfNotRegular}).
For each $\theta\in(0,1)$ and for each $\ell\in\left(\left[A_{0},A_{1}\right]_{\theta}\right)^{*}$
there exists a sequence $\left\{ y_{n}\right\} _{n\in\mathbb{N}}$
of elements in $X_{0}^{\#}\cap X_{1}^{\#}$ such that $\left\Vert y_{n}\right\Vert _{\left[X_{0}^{\#},X_{1}^{\#}\right]_{\theta}}=\left\Vert \ell\right\Vert _{\left(\left[A_{0},A_{1}\right]_{\theta}\right)^{*}}$
and 
\begin{equation}
\ell\left(a\right)=\lim_{n\to\infty}\left\langle a,y_{n}\right\rangle _{\left[X_{0},X_{1}\right]_{\theta}}\mbox{ for all }a\in\left[A_{0},A_{1}\right]_{\theta}.\label{eq:wibfy}
\end{equation}

\end{lem}

\noindent\textit{Proof .} Let $\ell$ be an arbitrary element of
$\left(\left[A_{0},A_{1}\right]_{\theta}\right)^{*}$ which we can
of course assume to be non-zero. We apply Lemma \ref{lem:3First}
to obtain the functions $h_{r}$ with the properties specified there.
Since $X_{0}^{\#}\cap X_{1}^{\#}$ is dense in $\left[X_{0}^{\#},X_{1}^{\#}\right]_{\theta}$,
for each $n\in\mathbb{N}$, there exists an element $v_{n}\in X_{0}^{\#}\cap X_{1}^{\#}$
such that $\left\Vert v_{n}-h_{1/n}\left(\theta\right)\right\Vert _{\left[X_{0}^{\#},X_{1}^{\#}\right]_{\theta}}\le1/n$.
Then, for each $a\in\left[A_{0},A_{1}\right]_{\theta}$, in view of
(\ref{eq:frtu}), 
\begin{eqnarray*}
\left|\left\langle a,v_{n}\right\rangle _{\left[X_{0},X_{1}\right]_{\theta}}-\left\langle a,h_{1/n}(\theta)\right\rangle _{\left[X_{0},X_{1}\right]_{\theta}}\right| & \le & \left\Vert a\right\Vert _{\left[X_{0},X_{1}\right]_{\theta}}\cdot\frac{1}{n}\,.
\end{eqnarray*}
This, together with (\ref{eq:oaists}), gives us that 
\begin{equation}
\ell(a)=\lim_{n\to\infty}\left\langle a,h_{1/n}(\theta)\right\rangle =\lim_{n\to\infty}\left\langle a,v_{n}\right\rangle .\label{eq:bztopa}
\end{equation}
We shall complete the proof by showing that 
\begin{equation}
\lim_{n\to\infty}\left\Vert v_{n}\right\Vert _{\left[X_{0}^{\#},X_{1}^{\#}\right]_{\theta}}=\left\Vert \ell\right\Vert _{\left(\left[A_{0},A_{1}\right]_{\theta}\right)^{*}},\label{eq:stbd}
\end{equation}
since this will imply that the sequence $\left\{ y_{n}\right\} _{n\in\mathbb{N}}$,
obtained by setting 
\[
y_{n}=\frac{\left\Vert \ell\right\Vert _{\left(\left[A_{0},A_{1}\right]_{\theta}\right)^{*}}}{\left\Vert v_{n}\right\Vert _{\left[X_{0}^{\#},X_{1}^{\#}\right]_{\theta}}}v_{n}
\]
has all the required properties. 

In view of (\ref{eq:bztopa}) and (\ref{eq:frtu}), we have that $\left|\ell(a)\right|\le\liminf_{n\to\infty}\left\Vert a\right\Vert _{\left[X_{0},X_{1}\right]_{\theta}}\left\Vert v_{n}\right\Vert _{\left[X_{0}^{\#},X_{1}^{\#}\right]_{\theta}}$
for each $a\in\left[A_{0},A_{1}\right]_{\theta}$. Taking the supremum
as $a$ ranges over the unit ball of $\left[A_{0},A_{1}\right]_{\theta}$,
we obtain that 
\begin{equation}
\left\Vert \ell\right\Vert _{\left(\left[A_{0},A_{1}\right]_{\theta}\right)^{*}}\le\liminf_{n\to\infty}\left\Vert v_{n}\right\Vert _{\left[X_{0}^{\#},X_{1}^{\#}\right]_{\theta}}.\label{eq:tac}
\end{equation}
Furthermore, by our choice of $v_{n}$ and by (\ref{eq:afyd}), 
\[
\left\Vert v_{n}\right\Vert _{\left[X_{0}^{\#},X_{1}^{\#}\right]_{\theta}}\le\left\Vert h_{1/n}(\theta)\right\Vert _{\left[X_{0}^{\#},X_{1}^{\#}\right]_{\theta}}+\frac{1}{n}\le e^{1/n}\left\Vert \ell\right\Vert _{\left(\left[A_{0},A_{1}\right]_{\theta}\right)^{*}}+\frac{1}{n}
\]
so that $\limsup_{n\to\infty}\left\Vert v_{n}\right\Vert _{\left[X_{0}^{\#},X_{1}^{\#}\right]_{\theta}}\le\left\Vert \ell\right\Vert _{\left(\left[A_{0},A_{1}\right]_{\theta}\right)^{*}}$.
This, combined with (\ref{eq:tac}), gives us (\ref{eq:stbd}) and
completes the proof of Lemma \ref{lem:Preparation}. $\qed$

Now we can easily provide the

\noindent \textit{Proof of Theorem \ref{thm:oldbaba}.} Given $a\in\left[X_{0},X_{1}\right]_{\theta}$,
which we can of course assume to be non-zero, we let $\ell$ be a
norm one element of $\left(\left[X_{0},X_{1}\right]_{\theta}\right)^{*}$
such that $\ell(a)=\left\Vert a\right\Vert _{\left[X_{0},X_{1}\right]_{\theta}}$.
We apply Lemma \ref{lem:Preparation} to obtain a sequence $\left\{ y_{n}\right\} _{n\in\mathbb{N}}$
in $X_{0}^{\#}\cap X_{1}^{\#}$ with $\left\Vert y_{n}\right\Vert _{\left[X_{0}^{\#},X_{1}^{\#}\right]_{\theta}}=1$
which satisfies (\ref{eq:wibfy}). Given $\varepsilon>0$, we can
now choose $b=y_{n}$ where, by virtue of (\ref{eq:wibfy}), $n$
can be chosen sufficiently large to ensure that $\left|\left\langle a,b\right\rangle _{\left[X_{0},X_{1}\right]_{\theta}}\right|\ge(1-\varepsilon)\ell(a)=\left(1-\varepsilon\right)\left\Vert a\right\Vert _{\left[X_{0},X_{1}\right]_{\theta}}$.
We also have $\left\langle a,b\right\rangle _{\left[X_{0},X_{1}\right]_{\theta}}=\left\langle a,b\right\rangle _{X_{0}+X_{1}}$
simply because $b$ is an element of $X_{0}^{\#}\cap X_{1}^{\#}=\left(X_{0}+X_{1}\right)^{\#}$
and (cf.~(\ref{eq:IntSp})) any sequence $\left\{ a_{n}\right\} _{n\in\mathbb{N}}$
in $X_{0}\cap X_{1}$ which converges to $a$ in the norm of $\left[X_{0},X_{1}\right]_{\theta}$
also converges in the norm of $X_{0}+X_{1}$. Thus we have established
all parts of (\ref{eq:faab}) except for the inequality on the right
side. But that is simply (\ref{eq:frtu}). $\qed$

\section{\protect Appendices}

Well, so far, there is only one appendix.

\subsection{\label{hilphi}Banach space valued analytic functions.}

\smallskip{}
 Let us here recall some basic definitions and properties of analytic
functions taking values in a Banach space. Our reference will be the
old treatise \cite{hp} of Hille and Phillips. We first need a preliminary
notion:
\begin{defn}
\label{def:DeterminingManifold}(Cf.\ \cite{hp} Definition 2.8.6,
p.\ 34.) Let $A$ be a Banach space and let $Y$ be a closed subspace
of $A^{*}$. Then $Y$ is said to be a \textit{``determining manifold''}
in $A^{*}$ if 
\begin{equation}
\left\Vert a\right\Vert _{A}=\sup\left\{ \left|y(a)\right|:y\in Y,\left\Vert y\right\Vert _{A^{*}}\le1\right\} \text{ for all }a\in A.\label{beawl}
\end{equation}

\end{defn}

\begin{defn}
\label{def:HolomVecValued}(Cf.\ \cite{hp} Definition 3.10.1, pp.\ 92--93.)
Let $A$ be a complex Banach 
space and let $\Omega\subset\mathbb{C}$ be a domain. An $A$ valued
function $f:\Omega\rightarrow A$ is said to be \textit{analytic}
(or \textit{holomorphic}) in $D$ if the scalar function $z\mapsto y(f(z))$
is analytic in $\Omega$ for each choice of the constant linear functional
$y$ in some determining manifold $Y$ in $A^{*}$. 
\end{defn}
Whenever we deal with this notion in this document, we will have the
option of choosing our determining manifold $Y$ to be either all
of $A^{*}$ or some predual of $A$ identified as a subspace of $A^{*}$
in the usual canonical way. (For example, when dealing with $X_{0}^{\#}+X_{1}^{\#}$
valued analytic functions, such as in the contexts of Remark \ref{rem:WeakerAnalyticity}
and the second part of the proof of Theorem \ref{babaeq}, we have
the convenient option of choosing $Y=X_{0}\cap X_{1}$, since $A=X_{0}^{\#}+X_{1}^{\#}=\left(X_{0}\cap X_{1}\right)^{\#}$.) 

Fortunately the property defined in Definition \ref{def:HolomVecValued}
turns out to be independent of the choice of the determining manifold
$Y$. This is because, perhaps surprisingly, the seemingly rather
weak condition imposed in that definition implies, and is therefore
equivalent, to a seemingly much stronger condition, namely that, for
each open disc $\{z:|z-z_{0}|<R\}$ contained in $\Omega$, there
exists a sequence $\left\{ a_{n}\right\} _{n\ge0}$ of elements of
$A$ such that the series $\sum_{n=0}^{\infty}(z-z_{0})^{n}a_{n}$
converges absolutely in $A$ norm to $f(z)$ at every point of that
disc. 

This follows from arguments (see \cite{hp} pp.\ 93--97) which first
show, using the uniform bounded principle, that any function $f:\Omega\to A$
satisfying the condition of Definition \ref{def:HolomVecValued} must
be continuous with respect to the norm topology of $A$. This continuity
enables one to use complex line integrals of $A$ valued continuous
functions to develop an obvious analogue of the theory of scalar valued
holomorphic functions. In particular, there are Cauchy integral formul\ae\
for $f$ and also for its derivatives of all orders, which are thus
revealed to exist and to also be $A$ valued holomorphic functions
on $\Omega$. In other words, for any function $f:\Omega\to A$ satisfying
the condition of Definition \ref{def:HolomVecValued}, there exist
functions $f^{(n)}:\Omega\to A$ which also satisfy that condition
and can be related to $f$ in various equivalent ways, not only by
the above-mentioned Cauchy integral formul\ae, but also, in particular,
by the formula 
\begin{equation}
y\left(f^{(n)}(\zeta)\right)=\frac{d^{n}}{d\zeta^{n}}\left(y(f(\zeta))\right)\mbox{ for each }\zeta\in\Omega\mbox{ and }y\in Y.\label{eq:wrnt}
\end{equation}

It is also sometimes useful to know that the convergence of difference
quotients of $f$ to the derivative of $f$ is also with respect to
the norm of $A$, i.e., that 
\begin{equation}
\lim_{h\to0}\left\Vert \frac{1}{h}\left(f(\zeta_{0}+h)-f(\zeta_{0})\right)-f^{\prime}(\zeta)\right\Vert _{A}=0\mbox{ for each }\zeta_{0}\in\Omega.\label{eq:NormCgce2Deriv}
\end{equation}
This can be immediately deduced from the fact, mentioned above, that
$f$ has a power series which converges in norm in a disc centred
at $\zeta_{0}$.
\begin{rem*}
Although this is not relevant for our discussion in these notes, let
us mention that if the function $f:\Omega\to A$ satisfies a weaker
version of the condition of Definition \ref{def:HolomVecValued} in
which the subspace $Y$ of $A^{*}$ is still required to satisfy (\ref{beawl})
but is no longer required to be closed, then the uniform boundedness
principle is not available for the above argument. In such a setting
it is necessary to impose other conditions to ensure that $f$ has
a power series representation etc. See for example \cite{ArendtWNikolskiN2000}
for more details of such matters.
\end{rem*}
\smallskip{}

\end{document}